\documentclass{article}
\usepackage{graphicx}

\usepackage{amssymb}
\newtheorem{theorem}{THEOREM}[section]
\newtheorem{proposition}{PROPOSITION}[section]
\newtheorem{lemma}{LEMMA}[section]
\newtheorem{corollary}{COROLLARY}[section]

\begin{document}
\title{The semiflow of a reaction diffusion equation with a singular potential}

\date{}
\author{Nikos. I. Karachalios, \\
{\it Department of Mathematics},\\
{\it University of the Aegean},\\
{\it Karlovassi, 83200 Samos, Greece}\\
\\
Nikolaos B. Zographopoulos,\\
{\it Department of Sciences}, \\
{\it Division of Mathematics}, \\
{\it Technical University of Crete}, \\
{\it 73100 Chania, Greece}
\\
}
\maketitle
\pagestyle{myheadings} \thispagestyle{plain} \markboth{N. B.
Zographopoulos}{Reaction diffusion equations with a singular
potential}
\maketitle
\begin{abstract}
We study the semiflow  $\mathcal{S}(t)$ defined by a semilinear
parabolic equation with a singular square potential
$V(x)=\frac{\mu}{|x|^2}$. It is known that the Hardy-Poincar\'{e}
inequality and its improved versions, have a prominent role on the
definition of the natural phase space.  Our study concerns the
case $0<\mu\leq\mu^*$, where $\mu^*$ is the optimal constant for
the Hardy-Poincar\'{e} inequality. On a bounded domain of
$\mathbb{R}^N$, we justify the global bifurcation of nontrivial
equilibrium solutions for a reaction term $f(s)=\lambda
s-|s|^{2\gamma}s$, with $\lambda$ as a bifurcation parameter. The
global bifurcation result is used to show that any solution
$\phi(t)=\mathcal{S}(t)\phi_0$, initiating form initial data
$\phi_0\geq 0$ ($\phi_0\leq 0$), $\phi_0\not\equiv 0$, tends to
the unique nonnegative (nonpositive) equilibrium.
\end{abstract}
\emph{Keywords:\ 
} \vspace{0.1cm} 
\section{Introduction}
Fundamental issues of the linear heat equation with a singular potential
\begin{eqnarray}
\label{linP}
\partial_t\phi-\Delta\phi&-&\frac{\mu}{|x|^2}\phi=0,\;\;x\in\Omega,\;t>0,\nonumber\\
\phi(x,0)&=&\phi_0(x),\;\;x\in\Omega,\\
\phi|_{\partial\Omega}&=&0,\;t>0,\nonumber
\end{eqnarray}
where $\Omega$ is in general an open set of $\mathbb{R}^N$, have been  analyzed in the works \cite{BarGo84, Cabre99, vz00}. The behavior of the solutions depends heavily on the critical value of the parameter $\mu$ (denoted by $\mu^*$),
which is the best constant of the Hardy's inequality.
The first fundamental result was that of \cite{BarGo84} for the Cauchy-Dirichlet problem in an open set of $\mathbb{R}^N$:
If $\phi_0(x)\geq 0$, $\phi_0(x)\neq 0$, there exists a global solution if $0<\mu\leq\mu_*$. Not even a local solution exists if $\mu>\mu^*$ (complete instantaneous blowup). The importance of Hardy's inequality for the result of \cite{BarGo84} was shown in \cite{Cabre99}. Further fundamental results ranging from the removal of the sign condition on the initial data, the uniqueness of solutions in an appropriate functional space, the possible decay of solutions and its rate when $\mu<\mu^*$, to the description of the behavior of solutions at the critical value $\mu^*$ as well as analysis of the Cauchy problem, have been addressed in \cite{vz00}. The legitimate analysis at the transition $\mu=\mu^*$ and beyond, for the aforementioned questions is related to the Hardy inequality  and to its improved versions \cite{vz00}, in bounded as well in unbounded domains.

Concerning the bounded domain case, the situation regarding the behavior of solutions of (\ref{linP}), can be described in summary as follows: When $0<\mu<\mu^*$ the sign-condition on the initial data can be removed and for any $\phi_0\in L^2(\Omega)$ there exist a unique, global in time solution $\phi\in C([0,\infty);L^2(\Omega))\cap L^2([0,\infty);H^1_0(\Omega))$ which decays at exponential rate. In the critical transition value $\mu=\mu^*$  solutions which still exist globally in $L^2(\Omega)$, {\em blow up instantaneously in $H^1_0(\Omega)$ but exist globally in the generalized Sobolev space $H_{\mu}^1 (\Omega)$}. The Hilbert space $H_{\mu}^1 (\Omega)$ is
defined for any fixed $0 < \mu \leq \mu^*$,  as the completion of
the $C_{0}^{\infty} (\Omega)$ functions under the norm
$||\phi||^{2}_{\mu} = \int_{\Omega} |\nabla \phi|^2\, dx -
\mu \int_{\Omega} \frac{\phi^2}{|x|^2}\, dx$. The $H^1_{\mu}(\Omega)$-solution, $\phi(x,t)\sim O(e^{-\lambda t})$ with the rate $\lambda>0$ explicitly given. This is only a basic framework, since there are important and deep phenomena in the range $0<\mu\leq\mu^*$ (e.g. singular behavior at the origin with prescribed rate, even for the solutions $\phi\geq 0$ associated with good initial data and even non uniqueness of nonnegative distributional solutions). When $\mu>\mu^*$, there exist initial data of oscillating type for which the solution exists globally in time. Extensions of the results of the bounded domain when $0<\mu\leq\mu^*$ (but with major differences e.g on the rate of decay) have been made on appropriate weighted spaces based on weighted improvements of the Hardy's inequality.

Strongly motivated by the results of \cite{vz00}, for (\ref{linP}) on the bounded domain case, we shall discuss the dynamics of a semilinear analogue of (\ref{linP})
\begin{eqnarray} \label{eq1.0a}
\partial_t\phi - \Delta \phi - \frac{\mu}{|x|^2}\phi&=& \lambda\, \phi
- |\phi|^{2\gamma}\,
\phi, \;\;x\in\Omega,\;\;t>0, \nonumber\\
\phi(x,0)&=& \phi_0(x),\;\;x\in\Omega,\\
\phi|_{\partial\Omega}&=&0,\;t>0\nonumber.
\end{eqnarray}
with our attention restricted in this work, up to the {\em critical case} $\mu=\mu^*$.

We start with the analysis of the set of equilibrium solutions of (\ref{eq1.0a}).
The equilibrium solutions in this case satisfy
the semilinear elliptic equation
\begin{eqnarray}
\label{eq1.1}
- \Delta u - \frac{\mu}{|x|^2}u &=& \lambda u - |u|^{2\gamma}u, \\
u|_{\partial \Omega} &=& 0, \nonumber
\end{eqnarray}
The results of Section 2, concern the bifurcation of equilibrium solutions with respect to the parameter $\lambda\in\mathbb{R}$. Considering this type of nonlinear term with $\lambda$ as a varying parameter, is of importance, having in mind the Ginzburg-Landau nonlinearity. Hardy's inequality implies  for the \emph{subcritical case} $0<\mu<\mu^*$, the equivalence $H^1_0(\Omega)\equiv H^1_{\mu}(\Omega)$. In this case, the operator $\mathcal{L}=-\Delta-\frac{\mu}{|x|^2}$ defines an unbounded self-adjoint operator in $L^2(\Omega)$ with compact inverse. Thus, a global branch of nonnegative solutions of (\ref{eq1.1}) bifurcating from the trivial solution at $(\lambda_1, 0)$, where $\lambda_1$ is the positive principal eigenvalue of the linear eigenvalue problem
\begin{eqnarray}
\label{eq1.2}
- \Delta u - \frac{\mu}{|x|^2}u &=& \lambda u, \\
u|_{\partial \Omega} &=& 0, \nonumber
\end{eqnarray}
is naturally expected.

The analysis carried out in \cite{vz00} for the critical case $\mu = \mu^*$, suggests that we cannot
expect $H^{1}_{0}$-solutions for the eigenvalue problem (\ref{eq1.1}). Instead, the main result of Section 2, is stated in the following

\begin{theorem} \label{globalbif}
Let\ $\Omega \subset \mathbb{R}^N$,\ $N \geq 3$,\ be a bounded
domain. Assume that $0 < \mu\leq\mu^*$, and that
\begin{eqnarray}
\label{cruc1} 0<\gamma\leq
\frac{Nq-2N+2q}{2(N-q)}:=\gamma_*,\;\;\mbox{for
any}\;\;\frac{2N}{N+2}<q<2.
\end{eqnarray}
Then, the principal eigenvalue $\lambda_{1,\mu}$\
of (\ref{eq1.2}) considered in $H^{1}_{\mu}(\Omega)$, is a bifurcating point of the problem
(\ref{eq1.1}) (in the sense of Rabinowitz) and
$C_{\lambda_{1,\mu}}$is a global branch of nonnegative $H^{1}_{\mu}(\Omega)$- solutions of (\ref{eq1.1}).
\end{theorem}
For comparison results and properties of the linear eigenvalue problem (\ref{eq1.2}), we will refer to \cite{davdup03}.

The global bifurcation results of Section 2 are of a twofold meaning. On the one hand, they establish the existence of a global branch $C_{\lambda_{1,\mu}}$, of nonnegative solutions in the critical value $\mu=\mu^*$. The global branch has the properties proved in
\setcounter{proposition}{1}
\begin{proposition}
\label{propbend}
Let\ $\Omega \subset \mathbb{R}^N$,\ $N \geq 3$,\ be a bounded
domain. Assume that $0 < \mu \leq \mu^*$ and that (\ref{cruc1}) holds. Then
(i)\ The global branch $C_{\lambda_{1,\mu}}$ bends to the
right of $\lambda_{1,\mu}$ (supercritical bifurcation) and it is
bounded for $\lambda$ bounded.

(ii)\ Every solution $u \in C_{\lambda_{1,\mu}}$ is the unique
nonnegative solution for the problem (\ref{eq1.1}).
\end{proposition}
On the other hand, it is well known that qualitative properties of a dynamical system may not depend continuously on the variation of the parameters. A first question of this nature can be addressed, regarding the behavior of global branches of nonnegative solutions possessed by (\ref{eq1.1}) in domains not containing the origin, $\Omega_r=\Omega\setminus B_r(0)$. We may consider $r>0$ as a parameter: Rabinowitz's theorem, can be applied to prove the existence of global branches $C_{\lambda_{1,\mu,r}}$ in $\mathbb{R}\times H_{0}^1(\Omega_r)$, for any $0<\mu\leq\mu^*$.  How these branches behave as $r\rightarrow 0$? A simple but careful analysis on the asymptotics of the eigenpairs $(\lambda, u_{\lambda,r})\in C_{\lambda_{1,\mu,r}}$ as $r\rightarrow 0$ (on the account of the properties of the space $H_{\mu}(\Omega)$ and the regularity results in $H_{\mu}(\Omega_r$), combined with the Whyburn's Theorem is used to prove
\setcounter{theorem}{2}
\begin{theorem} \label{approx}
Let\ $\Omega \subset \mathbb{R}^N$,\ $N \geq 3$,\ be a bounded
domain. Assume that $0 < \mu \leq \mu^*$ and that (\ref{cruc1}) holds. Then $C_{\lambda_{1,\mu,r}} \to C_{\lambda_{1,\mu}}$ in
$\mathbb{R}\times H_{\mu} (\Omega)$, as $r \downarrow 0$.
\end{theorem}

\begin{figure}
\begin{center}
    \begin{tabular}{cc}
    \includegraphics[scale=0.45]{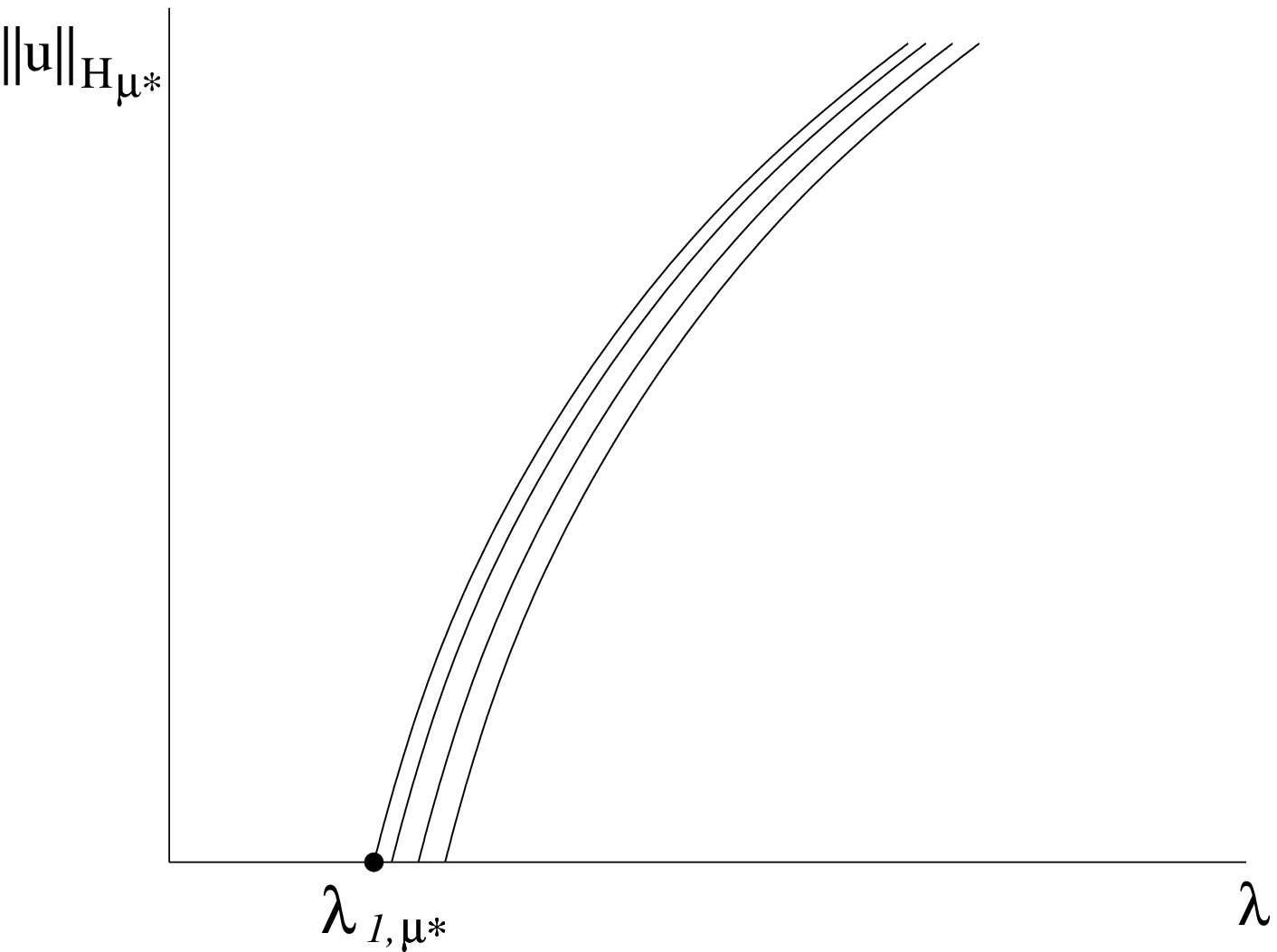} &
    \includegraphics[scale=0.45]{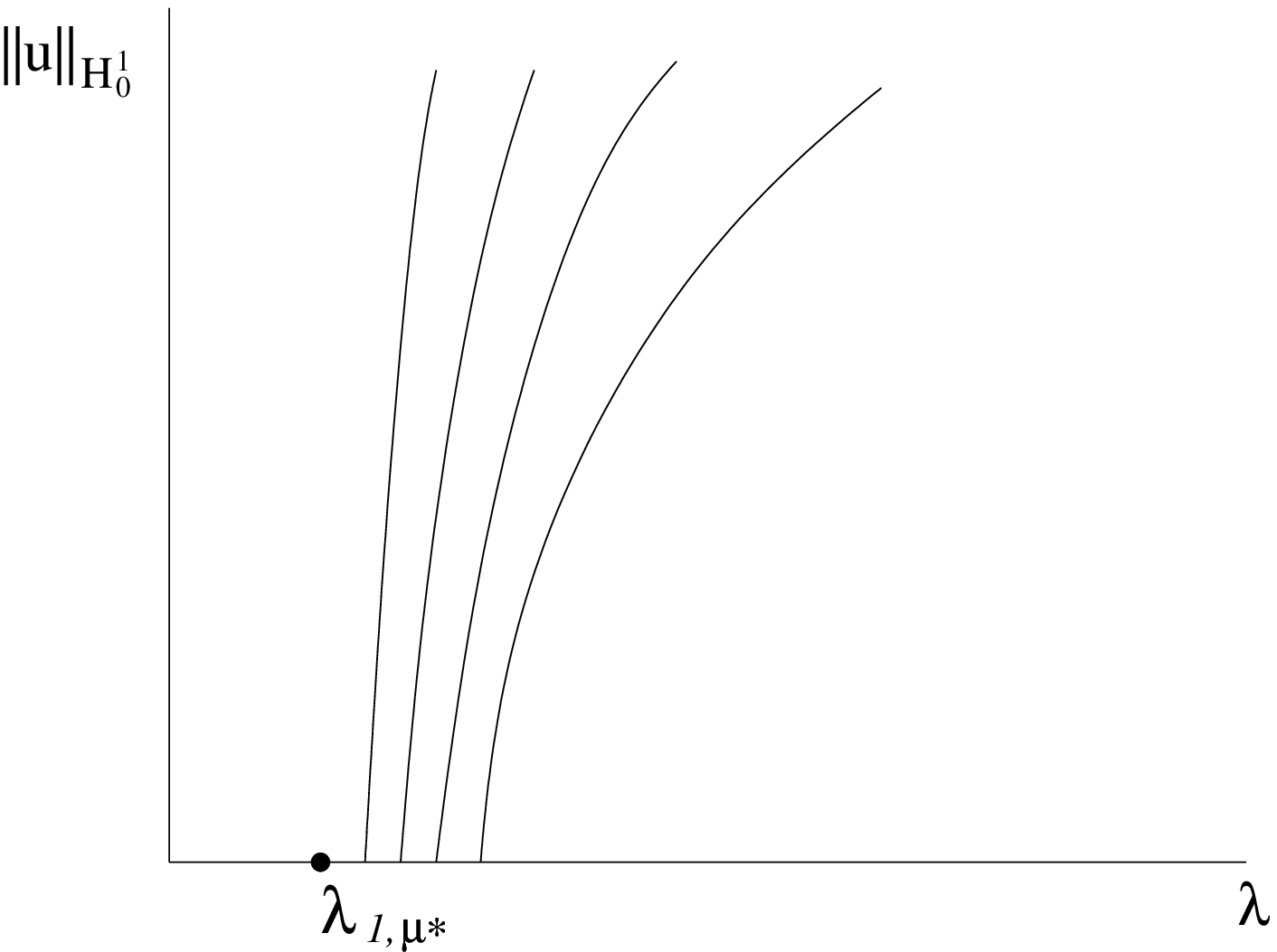}
    \end{tabular}
    \caption{(a) Possible bifurcation diagrams in (a) $H_{\mu^*} (\Omega)$  and (b) in $H^1_0(\Omega)$.}
    \label{fig1}
\end{center}
\end{figure}

It seems even  more interesting  to discuss {\em how the branches $C_{\lambda_{1,\mu}}$ for $\mu<\mu^*$ behave as $\mu\uparrow\mu^*$}.  Regarding the behavior of the global branch $C_{\lambda_{1,\mu}}$ as the parameter $\mu$ varies to the transition value $\mu^*$,
the answer is given in the following theorem, showing that the situation in $H^1_0(\Omega)$ and in $H^1_{\mu}(\Omega)$ is qualitatively totally different (figure \ref{fig1} demonstrates a possible configuration).
\setcounter{theorem}{3}
\begin{theorem} \label{convergth} Let $\Omega\subset\mathbb{R}^N$ be a bounded domain. We assume that
\begin{eqnarray}
\label{cruc1aN=3}
0<\gamma\leq\frac{5q-6}{2(3-q)}:=\gamma_*,\;\;\mbox{for any}\;\;\frac{3}{2}<q<2,\;\;\mbox{if}\;\;N=3,
\end{eqnarray}
and when $N\geq 4$, we assume condition (\ref{cruc1}).

A. Let $\mu_n \uparrow \mu^*$, as $n \to \infty$. Assume that
$(\lambda_n, u_n) \in C_{\lambda_{1,\mu_n}}$, be such that
$\lambda_n$ is bounded, i.e. $|\lambda_n|<L$. Then, $u_n$ must be
bounded too, in $H_{\mu^*} (\Omega)$. Moreover, $(\lambda_n, u_n)
\to (\lambda_*, u_*)$ in $\mathbb{R} \times H_{\mu^*} (\Omega)$,
with $(\lambda_*, u_*) \in C_{\lambda_{1,\mu^*}}$.
\newline
B. Let $\mu_n \uparrow \mu^*$, as $n \to \infty$. Assume that
$(\lambda_n, u_n) \in C_{\lambda_{1,\mu_n}}$, be such that
$\lambda_n \to \lambda_{1,\mu^*}$. Then, $u_n$ must be unbounded
in $H_0^1 (\Omega)$.
\end{theorem}
Observe that the condition (\ref{cruc1}) is slightly modified, distinguishing between the cases $N=3$ and $N\geq 4$. 

In Section 3, and in the spirit of our recent work \cite{kz06}, we shall use Theorem \ref{globalbif} to discuss the stability properties of equilibria and the asymptotic behavior of solutions of (\ref{eq1.0a}). We discuss first stability by linearization: Using the improved Hardy inequality of \cite{vz00} and its consequences, we consider appropriate Garding forms to prove the asymptotic stability of the trivial equilibrium when $\lambda\leq\lambda_{1,\mu}$ and the asymptotic stability of the unique nonnegative equilibrium when $\lambda >\lambda_{1,\mu}$, for $0<\mu\leq\mu^*$.  However the setting of \cite{vz00}, enables for a stronger result:  Following closely the semiflow theory  \cite{Ball2,jhale88,RTem88}, we define  a gradient semiflow in $H^1_{\mu}(\Omega)$, for any $0<\mu\leq\mu^*$.  This is one of the basic results proved in Section 3, stated in
\setcounter{proposition}{4}
\begin{proposition} \label{main1}
Let\ $\Omega \subset \mathbb{R}^N$,\ $N \geq 3$,\ be a bounded
domain, $0<\mu\leq\mu^*$ and condition (\ref{cruc1}) be fulfilled.
The semiflow (\ref{defsemiflow}),  possesses a global attractor
${\mathcal A}$ in $H_{\mu}(\Omega)$. Let $\mathcal{E}$ denote the
bounded set of equilibrium points of $\mathcal{S}(t)$. For each
complete orbit $\phi$ lying in $\mathcal{A}$, the limit sets
$\alpha(\phi)$ and $\omega(\phi)$ are connected subsets of $\mathcal{E}$ on which the
Lyapunov functional $\mathcal{J}$ associated to $\mathcal{S}(t)$,
is constant. If $\mathcal{E}$ is totally disconnected (in
particular if $\mathcal{E}$ is countable), the limit
$$z_{-}=\lim_{t\rightarrow-\infty}\phi(t),\;\;z_{+}=\lim_{t\rightarrow+\infty}\phi (t),$$
exist and are equilibrium points. Furthermore, any solution of
(\ref{eq1.0a}), tends to an equilibrium point as
$t\rightarrow\infty$.
\end{proposition}
Armed with the fact, that the limit set $\omega(\phi)$ for each positive orbit $\phi$ lying in the global attractor $\mathcal{A}$, is a connected subset of the bounded set $\mathcal{E}$ of the equilibrium solutions, the global bifurcation result of Theorem \ref{globalbif} will be crucial: It actually shows that $\mathcal{E}=\{0\}$ when $\lambda\leq\lambda_{1,\mu}$, and is totally disconnected  when $\lambda>\lambda_{1,\mu}$, $\mathcal{E}=\{u_{-},0,u\}$, $u_{-}=-u$ in $H_{\mu}(\Omega)$, for any $0<\mu\leq\mu^*$. The trivial solution is unstable when $\lambda>\lambda_{1,\mu}$, thus the limit set $\omega(\phi_0)$ for every $\phi_0\in H_{\mu}(\Omega)$ of definite sign, is described for all $0<\mu\leq\mu^*$ by
\setcounter{theorem}{5}
\begin{theorem}  \label{tsibadyoball} Let\ $\Omega \subset \mathbb{R}^N$,\ $N \geq 3$,\ be a bounded
domain. Assume that $0 < \mu \leq \mu^*$ and that (\ref{cruc1}) is fulfilled.  Let $\phi_0\in H_{\mu}(\Omega)$, $\phi_0\not\equiv 0$.  If $\lambda\leq\lambda_{1,\mu}$, then $\mathcal{A}=\{0\}$. If $\lambda>\lambda_{1,\mu}$, then $\omega({\phi_0})=\{u\}$ when $\phi_0\geq 0$ and $\omega(\phi_0)=\{u_{-}\}$ when $\phi_0\leq 0$.
\end{theorem}

\begin{figure}
\begin{center}
    \begin{tabular}{cc}
    \includegraphics[scale=0.7]{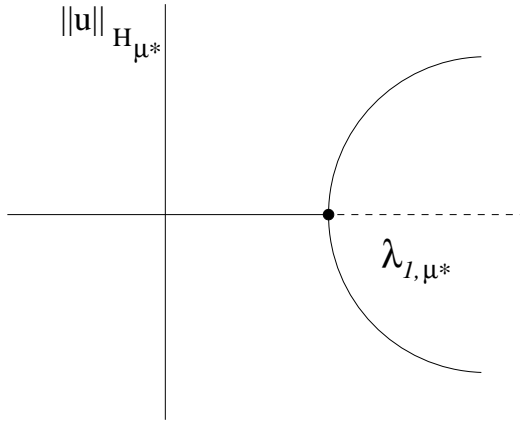}
    \end{tabular}
    \caption{Supercritical pitchfork bifurcation for the semiflow defined by (\ref{eq1.0a}) in $H_{\mu^*} (\Omega)$.}
    \label{fig2}
\end{center}
\end{figure}

The above result is a rigorous verification that (\ref{eq1.0a}) which undergoes a pitchfork bifurcation of supercritical type for any $\mu<\mu^*$ in $H^1_0(\Omega)$ preserves this behavior up to the transition $\mu=\mu^*$ in the $H^{1}_{\mu^*}(\Omega)$-phase space (see figure \ref{fig2}). We remark that in the case $\lambda>\lambda_{1,\mu}$ Proposition \ref{main1}, clearly implies that for {\em any} $\phi_0\not\equiv{0}$, {\em any} solution $\phi(t)=\mathcal{S}(t)\phi_0$ converges to one of the equilibrium solutions $u$ or $u_{-}$, possibly through an heteroclinic orbit connecting them.

However, Theorem \ref{convergth} B. combined with Theorem \ref{tsibadyoball} indicate for the ``explosive'' behavior of the attractor $\mathcal{A}$ in $H^1_{0}(\Omega)$ when $\mu\to\mu^*$. Theorem \ref{tsibadyoball} could also be viewed as the analogue of \cite[Theorem 4.1, pg. 123]{vz00} for (\ref{eq1.0a}), with the exponential decay, replaced by the convergence to the unique nonnegative or the unique nonpositive equilibrium, for any $\lambda>\lambda_{1,\mu}$, according to the sign of the initial data $\phi_0$.

At this point, we also remark \cite{chav06} for bifurcation results on $H^{1}_0(\Omega)$ with $\mu$ as a bifurcation parameter, regarding the semilinear elliptic problem
\begin{eqnarray*}
- \Delta u - \frac{\mu}{|x|^2}u = u^q,\;\;u>0,\;\;u|_{\partial \Omega} =0.
\end{eqnarray*}
For bifurcation results on the degenerate elliptic problem
\begin{eqnarray*}
- |x|^2\Delta u  = \lambda f(u),\;\;u>0,\;\;u|_{\partial \Omega} =0,
\end{eqnarray*}
related to the Hardy inequality, we refer to \cite{eg00}. We also point out \cite{barbabrown}, on recent bifurcation results for the elliptic problem
\begin{eqnarray*}
-\Delta u  = \lambda m(x)u+b(x)u^{\gamma},\;\;\frac{\partial u}{\partial n}|_{\partial \Omega} =0,
\end{eqnarray*}
where the functions $m,b:\overline{\Omega}\to\mathbb{R}$ are this time, smooth functions but of changing sign.
For a brief reference to existing results on the issue of convergence of solutions of global solutions of evolution equations to steady states, we refer to \cite{BJP02} (see also \cite[pg. 366]{kz06}). For improvements related to  second order Hardy-type inequalities, we refer to the recent work  \cite{TZ07}.
\section{Global bifurcation of equilibrium solutions}
\setcounter{equation}{0}
This section is devoted to the proof of the existence of  bifurcation branches for the equilibrium solutions of (\ref{eq1.0a}) given by the semilinear elliptic equation (\ref{eq1.1}).
Here $\Omega$ will be an open bounded and connected subset of
$\mathbb{R}^N$, $N \geq 3$ including the origin. We shall
assume that $0< \mu \leq \mu^*$, where
\[
\mu^* := \biggr ( \frac{N-2}{2} \biggr )^2,
\]
is the best constant of Hardy's inequality
\begin{equation} \label{hardyineq}
\int_{\Omega} |\nabla u|^2\, dx > \biggr ( \frac{N-2}{2} \biggr
)^2 \int_{\Omega} \frac{u^2}{|x|^2}\, dx.
\end{equation}
In subsection \ref{prelim} we recall the
basic properties of the delicate functional framework developed in
\cite[Section 4, pg. 121-123]{vz00},  and we present some auxiliary results regarding the nonlinear maps defined in this setting. Subsection \ref{THEOREMA2} refers to the proof of Theorem \ref{globalbif}, while subsection \ref{THEOREMA2} is devoted to the approximation of the global branch by the associated branches of systems considered in domains not containing the origin. In subsection \ref{THEOREMA3} we discuss the proof of Theorem \ref{convergth}.
\subsection{Basic properties of the phase space.}
\label{prelim}
It well  known that the constant $\mu^*$ is optimal  and it is
not attained in $H^{1}_{0} (\Omega)$. In  \cite{bv97} it was given the following improved version of (\ref{hardyineq})
\begin{equation} \label{HaPoi1}
\int_{\Omega} |\nabla u|^2\, dx \geq \biggr ( \frac{N-2}{2} \biggr
)^2 \int_{\Omega} \frac{u^2}{|x|^2}\, dx + \lambda_{\Omega}
\int_{\Omega} u^2\, dx,
\end{equation}
where $\lambda_{\Omega} = z_0^2\, \omega_N^{\frac{2}{N}}\,
|\Omega|^{-\frac{2}{N}}$, where $\omega_N$  and $|\Omega|$ denote
the volume of the unit ball and $\Omega$ respectively, and $z_0 =
2.4048 \ldots$ denotes the first zero of the Bessel function
$J_0(z)$. This constant is optimal when $\Omega$ is a ball, but it
is also not achieved in $H_{0}^{1} (\Omega)$.  In
\cite{ft02} was proved that inequality (\ref{hardyineq}) admits an
infinite series of correction terms.

The analysis of  \cite{vz00},  recovered that the natural
phase space for the study of linear equation (\ref{linP})
system (\ref{eq1.0a}) is the Hilbert space $H_{\mu} (\Omega)$,
defined for any fixed $0 < \mu \leq \mu^*$,  as the completion of
the $C_{0}^{\infty} (\Omega)$ functions under the norm
\begin{eqnarray}
\label{ornorm}
||\phi||^{2}_{\mu} = \int_{\Omega} |\nabla \phi|^2\, dx -
\mu \int_{\Omega} \frac{\phi^2}{|x|^2}\, dx,
\end{eqnarray}
and endowed with the scalar product
\[
\left(\phi,\psi\right)_{\mu} = \int_{\Omega} \nabla \phi
\nabla \psi\, dx - \mu \int_{\Omega} \frac{\phi\, \psi}{|x|^2}\,
dx.
\]
Consequently, this is also the case for the semilinear analogue (\ref{eq1.0a}).
Friedrich's extension theory is applicable due to the  inequality (\ref{HaPoi1}): is the main ingredient which can be used to consider the operator $\mathcal{L}=-\Delta-V(x)$ as a positive
and self adjoint operator with domain of definition
\begin{eqnarray}
\label{HaPoi2} D(\mathcal{L})=\left\{\phi\in H_{\mu}
(\Omega):\mathcal{L}\phi\in L^{2}(\Omega)\right\}.
\end{eqnarray}
The {\em improved Hardy-Poincar\'{e} inequalities}
\begin{eqnarray}
\label{HaPoiim1}
\int_{\Omega}\left[|\nabla\phi-\mu^*\frac{\phi^2}{|x|^2}\right]dx&\geq&
C(q,\Omega)||\phi||_{W^{1,q}(\Omega)}^2,\;\;1\leq q<2,\\
\label{HaPoiim2}
\int_{\Omega}\left[|\nabla\phi-\mu^*\frac{\phi^2}{|x|^2}\right]dx&\geq&
C(s,r, \Omega)||\phi||_{W^{s,r}(\Omega)}^2,\;\;0\leq s<1,\;\;1\leq
r<r_{*}=\frac{2N}{N-2(1-s)},
\end{eqnarray}
for all $\phi\in C^{\infty}_0(\Omega)$, imply the
continuous embeddings,
\begin{eqnarray}
\label{HaPoiim3} H_{\mu}(\Omega) \hookrightarrow W_{0}^{1,q}
(\Omega),\;\;\; H_{\mu}(\Omega) \hookrightarrow H_{0}^{s}
(\Omega),\;\;1 \leq q<2,\;\;0\leq s<1.
\end{eqnarray}
if $1 \leq q<2$ and $0 \leq s<1$. Furthermore, since $W_{0}^{1,q}
(\Omega)$ is compactly embedded in $H_{0}^{s}$ for suitable
$q=q(s)$ close enough to $2$, and $H_{0}^{s} (\Omega)$ is
compactly embedded in $L^2(\Omega)$, we infer the compact
embeddings
\begin{eqnarray}
\label{HaPoiim4} H_{\mu}(\Omega) \hookrightarrow\hookrightarrow
L^2(\Omega),\;\; H_{\mu}(\Omega) \hookrightarrow\hookrightarrow
H_{0}^{s}(\Omega),\;\;0\leq s<1.
\end{eqnarray}

In the subcritical case $0 < \mu < \mu^*$  we have the following property
of $H_{\mu} (\Omega)$.
\begin{lemma} (\cite{vz00}) \label{equivspace}
Let $0<\mu<\mu^*$. Then $H_{\mu} (\Omega) \equiv H^{1}_{0}(\Omega)$.
\end{lemma}
{\bf Proof:} Clearly from (\ref{ornorm}),
\begin{equation}\label{eq2.2}
||u||^2_{\mu} \leq \int_{\Omega} |\nabla u|^2\, dx =
||u||^2_{H^{1}_{0}(\Omega)}.
\end{equation}
On the other hand, Hardy's inequality (\ref{hardyineq}), implies that
\begin{equation}\label{eq2.3}
||u||^2_{{\mu}} \geq \biggl [ 1 - \biggr ( \frac{N-2}{2}
\biggr )^{-2}\, \mu \biggr ]\, ||u||^2_{H^{1}_{0}(\Omega)}.
\end{equation}
Thus, from inequalities (\ref{eq2.2}) and (\ref{eq2.3}) we
conclude that

\begin{eqnarray*}
c\, ||u||^2_{H^{1}_{0}(\Omega)} \leq ||u||^2_{H_{\mu}(\Omega)} \leq
C\, ||u||^2_{H^{1}_{0}(\Omega)},
\end{eqnarray*}
for $c = 1 - \biggr ( \frac{N-2}{2} \biggr
)^{-2}\, \mu>0$ if $\mu<\mu^*$,  and $C=1$. \ $\blacksquare$
\newline
A remarkable property was shown in  \cite{vz00} concerning the
critical case $\mu=\mu^*$: $H_{\mu^*} (\Omega)$ is larger than
$H_{0}^1 (\Omega)$, since it contains singularities of the form $f
\sim |x|^{(N-2)/2}$, and it is smaller than $\cap_{q<2} W^{1,q}
(\Omega)$.

With the continuous embeddings (\ref{HaPoiim3}) at hand, we can handle the nonlinearity of (\ref{eq1.0a}).
\setcounter{lemma}{1}
\begin{lemma}
\label{aux1} Let condition (\ref{cruc1}) be satisfied and  assume that $\mu\leq\mu*$. The function
$g(s)=|s|^{2\gamma}s, s\in\mathbb{R}$, defines a sequentially
weakly continuous map $g:H_{\mu}(\Omega)\rightarrow L^2(\Omega)$.
Let $G(\phi):=\int_{0}^{\phi} g(s)ds$. The functional
$\mathcal{G}:H_{\mu}(\Omega)\rightarrow\mathbb{R}$ defined by
$\mathcal{G}(\phi)=\int_{\Omega}G(\phi)dx$, is $C^1$ and
sequentially weakly continuous.
\end{lemma}
{\bf Proof:} Starting by the standard
Sobolev embeddings,  we recall that
\begin{eqnarray}
\label{KZ0} W^{1,q}(\Omega)\hookrightarrow
L^{p}(\Omega)\;\;\mbox{for any}\;\;1\leq
p\leq\frac{qN}{N-q},\;\;q<N.
\end{eqnarray}
We consider the critical exponent
\begin{eqnarray}
\label{KZ1}
p^*:=\frac{qN}{N-q}\;\;\mbox{for any}\;\;1\leq q<2.
\end{eqnarray}
Thus, as an immediate consequence of the embedding
(\ref{HaPoiim3}) we infer that
\begin{eqnarray}
\label{KZ2} H_{\mu}(\Omega)\hookrightarrow
L^{p}(\Omega),\;\;\mbox{for any}\;\;1\leq p\leq p^*.
\end{eqnarray}
Using (\ref{KZ2}) it can be easily deduced that the functional $g$
is well defined, under the restriction (\ref{cruc1}). Furthermore,
it follows from (\ref{KZ2}), that $\mathcal{G}$ is well defined if
\begin{eqnarray*}
\label{dontworry} 0<\gamma\leq
\frac{Nq-2N+2q}{(N-q)}:=\gamma^*,\;\;\mbox{for
any}\;\;\frac{2N}{N+2}<q<2.
\end{eqnarray*}
noting that $\gamma_*<\gamma^*$.

That both functional are sequentially weakly continuous, can be
verified by using the compact embeddings (\ref{HaPoiim4}) and
repeating the arguments  of \cite[Lemma 3.3, pg. 38 \& Theorem
3.6, pg. 40]{Ball2}. We will check that $\mathcal{G}$ is a
$C^1$-functional, and its derivative is
\begin{eqnarray}
\label{weak5}
\mathcal{G}'(\phi)(z)=\left<g(\phi),z\right>,\;\;\mbox{for
every}\;\;\phi\in H_{\mu}(\Omega),\;\;z\in H_{\mu}^{-1}(\Omega).
\end{eqnarray}
We consider for $\phi,\psi\in H_{\mu}(\Omega) $, the quantity
\begin{eqnarray}
\;\;\;\;\;\;\;\;\frac{\mathcal{G}(\phi+s\psi)-\mathcal{G}(\phi)}{s}&=&\frac{1}{s}
\int_{\Omega}\int_{0}^{1}\frac{d}{d\theta}G(\phi+\theta
s\psi)d\theta dx
\nonumber\\
&=&\int_{\Omega}\int_{0}^{1}g(\phi+s\theta\psi)\psi d\theta dx.
\end{eqnarray}
We set  $\sigma=\frac{qN}{N(q-1)+q}$, $\sigma^{-1}+{p^*}^{-1}=1$, and we get
\begin{eqnarray}
\left|\int_{\Omega}g(\phi+\theta s\psi)\psi dx\right|\leq
c\left(\int_{\Omega}(|\phi|^{(2\gamma +1)}+|\psi|^{(2\gamma
+1)})^{\sigma}dx\right)^{\frac{1}{\sigma}}
\left(\int_{\Omega}|\psi|^{p^*}dx\right)^{\frac{1}{p^*}}.
\end{eqnarray}
To apply the continuous embedding (\ref{KZ2}) we need the
requirement $$(2\gamma +1)\sigma\leq p^*.$$ This requirement
produces the restriction (\ref{cruc1}). Letting $s\rightarrow 0$, and using the
dominated convergence theorem, we infer that $\mathcal{G}$ is differentiable
with the derivative (\ref{weak5}).

For the continuity, we consider a sequence
$\{\phi_n\}_{n\in\mathbb{N}}$ of $H_{\mu}(\Omega)$ such that
$\phi_n\rightarrow\phi$ in $H_{\mu}(\Omega)$ as $n\rightarrow
\infty$. We note first, that
\begin{eqnarray}
\label{weak6}
\left<\mathcal{G}'(\phi_n)-\mathcal{G}'(\phi),z\right>\leq
||g(\phi_n)-g(\phi)||_{L^\sigma}||z||_{L^{p^*}}.
\end{eqnarray}
Setting  then $p_1=\frac{p^*}{\sigma}$, the requirement for
$p_1>1$, produces again the restriction for $\frac{2N}{N+2}<q<2$.
Now for
$$p_2=\frac{N(q-1)+q}{N(q-1)-N+2q},\;\;p_2^{-1}+p_1^{-1}=1,$$
we get the inequality
\begin{eqnarray*}
||g(\phi_n)-g(\phi)||_{L^\sigma}^\sigma&\leq& c
\left(\int_{\Omega}(|\phi_n|^{2\gamma}+|\phi|^{2\gamma})^{\sigma
p_2}
dx\right)^{\frac{1}{p_2}}
\left(\int_{\Omega}|\phi_n-\phi|^{p^*}dx\right)^{\frac{1}{p_1}}.
\end{eqnarray*}
The embedding (\ref{KZ2}) is applicable if $2\gamma \sigma p_2<
p^*$, giving (\ref{cruc1}). Under this condition and as
\begin{eqnarray*}
\lim_{n\rightarrow\infty}\int_{\Omega}|\phi_n-\phi|^{p^*}dx=0,
\end{eqnarray*}
we conclude from (\ref{weak6}), the continuity of $\mathcal{G}'$.
\ $\blacksquare$

\subsection{Existence of a global branch of positive solutions for any $0<\mu\leq\mu^*$}
\label{THEOREMA1}
The existence of a global branch of nonnegative solutions will be proved via the classical Rabinowitz's theorem:
\setcounter{theorem}{2}
\begin{theorem} \label{rab}
Assume that $X$ is a Banach space with norm\ $||\cdot||$ and
consider\ $G(\lambda,\cdot)=\lambda L\cdot + H(\lambda,\cdot)$,\
where\ $L$\ is a compact linear map on\ $X$\ and\
$H(\lambda,\cdot)$\ is compact and satisfies
\begin{eqnarray}
\label{Order}
\lim_{||u|| \to 0} \frac{||H(\lambda,u)||}{||u||}=0.
\end{eqnarray}
If\ $\lambda$\ is a simple eigenvalue of\ $L$\ then the closure of
the set
\[
C=\{ (\lambda ,u) \in \mathbb{R} \times X : (\lambda,u)\;\;
\mbox{solves}\;\; u=G(\lambda,u),\; u \not\equiv 0 \},
\]
possesses a maximal continuum (i.e. connected branch) of
solutions,\ $C_\lambda$,\ such that\ $(\lambda,0) \in
C_{\lambda}$\ and\ $C_{\lambda}$\ either:

(i)\ meets infinity in\ $\mathbb{R} \times X$\ or,

(ii)\ meets\ $(\lambda^*,0)$,\ where\ $\lambda^* \ne \lambda$\ is
also an eigenvalue of\ $L$.
\end{theorem}
We will prove that there exists a global branch (i.e. the
second alternative of Theorem \ref{rab} cannot happen) of
solutions bifurcating from the principal eigenvalue
$\lambda_{1,\mu}$ of the problem (\ref{eq1.2}),
for any $\mu \leq \mu^*$.
\setcounter{lemma}{3}
\begin{lemma} \label{princeig}
Assume that $0 < \mu \leq \mu^*$. Problem (\ref{eq1.2}), admits a
positive principal eigenvalue\ $\lambda_{1,\mu}$, given by
\begin{equation}\label{varchareig}
\lambda_{1,\mu} = \inf_{
\begin{array}{c}
               \phi \in H_{\mu}(\Omega) \\
                \phi \not\equiv 0
               \end{array}}
\frac{\int_{\Omega} |\nabla \phi|^2\; dx - \mu \int_{\Omega}
\frac{\phi^2}{|x|^2}\, dx}{\int_{\Omega} |\phi|^2\; dx}.
\end{equation}
with the following properties:

(i) $\lambda_{1,\mu}$ is simple with a positive associated
eigenfunction\ $u_{1,\mu}$, which belongs at least to
$C^{1,\zeta}_{loc} (\Omega \backslash \{0\})$, for some\ $ \zeta
\in (0,1)$,

(ii) $\lambda_{1,\mu}$ is the only eigenvalue of (\ref{eq1.2})
with nonnegative associated eigenfunction.
\end{lemma}
{\bf Proof:} The existence and the variational characterization
(\ref{varchareig}) of the principal eigenvalue follows from the compactness of the embeddings (\ref{HaPoiim4}) implying that $\mathcal{L}=-\Delta-\frac{\mu}{|x|^2}$  for $0<\mu\leq\mu^*$, has an
orthonormal basis of eigenfunctions in $H_{\mu}(\Omega)$ with an
eigenvalue sequence
\begin{eqnarray}
\label{eigensec}
0<\lambda_1\leq\lambda_2\leq\cdots\leq\lambda_n\leq\cdots\rightarrow\infty,
\end{eqnarray}
(cf.
\cite[pg. 122]{vz00})
The regularity results (cf. \cite[Theorem 8.22]{giltru77})
imply that if\ $u$\ is a weak solution of the problem
(\ref{eq1.2}), then\ $u \in C^{2,\zeta}_{loc} (\Omega \backslash
\{0\})$, for some\ $ \zeta \in (0,1)$.\ The positivity of
$u_{1,\mu}$ follows from \cite[Lemma 2.2]{davdup03}-we also refer to the weak maximum principle of
\cite{brecab98}. The simplicity and the uniqueness up to positive
eigenfunctions of $\lambda_{1,\mu}$ can be verified, by using Picone's identity \cite{kz06}.\
$\blacksquare$ \vspace{0.2cm}

For some further properties of the principal eigenvalue and the
corresponding eigenfunction, we refer to \cite{davdup03}. We remark
\cite{ft02}, where the weighted space Hilbert space $W_0^{1,2}
(\Omega;\; |x|^{-(N-2)}$ was used, defined as the completion of
$C_{0}^{\infty}$-functions under the norm
\[
||u||_{W_0^{1,2} (\Omega;\, |x|^{-(N-2)})} = \int_{\Omega}
|x|^{-(N-2)} |\nabla w|^2\, dx + \int_{\Omega} |x|^{-(N-2)} w^2\,
dx
\]
and endowed with the inner product
\[
<u,v>_{W_0^{1,2} (\Omega;\, |x|^{-(N-2)})} = \int_{\Omega}
|x|^{-(N-2)} \nabla f \nabla g\, dx + \int_{\Omega}
|x|^{-(N-2)} f\, g\, dx.
\]
In \cite{ft02}, the space $W_0^{1,2}
(\Omega;\; |x|^{-(N-2)}$ was considered for the proof of the existence of principal eigenvalues for the
eigenvalue problem
\begin{eqnarray} \label{ft}
- \Delta u - \mu\, \frac{u}{|x|^2} &=& \lambda\, V(x)\, u \\
u|_{\partial \Omega} &=& 0, \nonumber
\end{eqnarray}
Furthermore, it was assumed that $V(x)\geq 0$, $V(x)\in L^p(\Omega)$, $p=N/2$. A comparison of the spaces $H_{\mu^*}(\Omega)$ and  $W_0^{1,2} (\Omega; |x|^{-(N-2)})$ implies that $u \in H_{\mu^*} (\Omega)$ if and only
if $|x|^{(N-2)/2} u \in W_0^{1,2} (\Omega; |x|^{-(N-2)})$.

Proceeding to the proof of the global bifurcation result, we discuss first the behavior of $\lambda_{1,\mu}$, $0<\mu<\mu^*$ as $\mu \uparrow
\mu^*$. Next lemma demonstrates the qualitative differences in $H^1_0(\Omega)$ for the solutions of the linear eigenvalue problem  (\ref{eq1.2}) as $\mu$ converges to the transition value $\mu^*$.
\setcounter{proposition}{4}
\begin{proposition}
\label{limit}
Let $\mu \uparrow \mu^*$. Then,

(i) $\lambda_{1,\mu}$ is a decreasing sequence, and there exists $\lambda_*>0$ such that
$\lambda_{1,\mu} \downarrow \lambda_*$.

(ii)\ The corresponding normalized eigenfunctions $u_{1,\mu}$ are
converging weakly to 0, in $H_{0}^{1} (\Omega)$.
\end{proposition}
{\bf Proof:} (i)\ Let $\mu_1 < \mu_2$. Then the variational characterization of the principal eigenvalue $\lambda_{1,\mu}$ (\ref{varchareig})
implies that $\lambda_{1,\mu_1} > \lambda_{1,\mu_2}$. Thus
$\lambda_{1,\mu}$ is decreasing. Applying next the improved Hardy's
inequality (\ref{HaPoi1}) we infer that $\lambda_{1,\mu}$ is bounded from below by $\lambda_{\Omega}$. Thus, there exists $\lambda_*>0$, such that $\lambda_{1,\mu}\downarrow\lambda_*$.\

(ii) The eigenfunctions $u_{1,\mu}$ should satisfy the weak formula
\begin{equation}\label{eq2.7}
\int_{\Omega} \nabla u_{1,\mu} \nabla \phi\, dx - \mu
\int_{\Omega} \frac{u_{1,\mu}\, \phi}{|x|^2}\, dx =
\lambda_{1,\mu} \int_{\Omega} u_{1,\mu}\, \phi\, dx,
\end{equation}
for any $\phi \in C_{0}^{\infty} (\Omega)$. We still denote by $u_{1,\mu}$ the sequence of normalized eigenfunctions, forming a bounded sequence in $H^1_0(\Omega)$.  We deduce that there exists some $u\in H^1_0(\Omega)$ such that up to a subsequence (not relabelled), $u_{1,\mu} \rightharpoonup u$ in
$H_{0}^{1} (\Omega)$ and $u_{1,\mu} \to u$ in $L^{q}(\Omega)$, for
any $1<q<\frac{2N}{N-2}$. For some fixed $\varepsilon>0$,
small enough and any $\phi \in C_{0}^{\infty} (\Omega)$, we have that
\[
\int_{\Omega} \frac{(u_{1,\mu} - u)\, \phi}{|x|^2}\, dx \leq
||\phi||_{L^{\infty}(\Omega)}\, \left ( \int_{\Omega} |u_{1,\mu} -
u|^{\frac{N-\varepsilon}{N-2-\varepsilon}} \right
)^{\frac{N-2-\varepsilon}{N-\varepsilon}}\, \left ( \int_{\Omega}
|x|^{-N+\varepsilon} \right )^{2/(N-\varepsilon)} \to 0,
\]
thus
\[
\int_{\Omega} \frac{u_{1,\mu}\, \phi}{|x|^2}\, dx \to
\int_{\Omega} \frac{u\, \phi}{|x|^2}\, dx,\;\;\; \mbox{as}\;\; \mu
\uparrow \mu^*.
\]
Let us now assume by contradiction that $u \not\equiv 0$. Passing to the limit in (\ref{eq2.7}), we
get that $u$ must satisfy
\[
\int_{\Omega} \nabla u \nabla \phi\, dx - \mu^*
\int_{\Omega} \frac{u\, \phi}{|x|^2}\, dx = \lambda_*
\int_{\Omega} u\, \phi\, dx,
\]
for any $\phi \in C_{0}^{\infty} (\Omega)$, or equivalently that $u$
must be a nontrivial solution of the problem
\begin{eqnarray}
\label{contra}
- \Delta u - \mu^* \frac{u}{|x|^2} = \lambda_*\, u,\;\;\; u \in
H_{0}^{1} (\Omega).
\end{eqnarray}
However, since $\mu^*$ is the optimal constant of (\ref{HaPoi1}) which is not achieved in $H^{1}_{0}(\Omega)$, (\ref{contra}) implies that $u\equiv 0$. $\blacksquare$
\vspace{0.2cm} \\
{\bf Proof of Theorem \ref{globalbif}:} For the justification of Theorem \ref{rab}, the improved Hardy's
inequality (\ref{HaPoi1}), will allow us to employ the method developed in \cite{bro94}: On the account of (\ref{varchareig}), we define a bilinear form  in
$C^{\infty}_0(\Omega)$ by
\begin{eqnarray} \label{Rabi1}
<u,v>_X = \int_{\Omega} \nabla u \nabla v \, dx - \mu\,
\int_{\Omega} \frac{u\, v}{|x|^2}\, dx - \frac{c}{2} \int_{\Omega}
u\, v\, dx, \;\; \mbox{for all}\;\; u, v \in C^{\infty}_0(\Omega),\;\;c = \lambda_{1,\mu}.
\end{eqnarray}
We define next the space $X$, as
the completion of $C^{\infty}_0(\Omega)$ with respect to the
norm induced by  (\ref{Rabi1}), $||u||_{X}^2= <u,u>_X$: Then due to the improved Hardy's
inequality (\ref{HaPoi1}), we deduce the equivalence of norms
\begin{eqnarray*}
\frac{1}{2} ||u||^2_{H_{\mu}(\Omega)} \leq ||u||_{X}^2 \leq
\frac{3}{2} ||u||^2_{H_{\mu}(\Omega)},\;\; \mbox{for all} \;\; u,
v \in C^{\infty}_0(\Omega),
\end{eqnarray*}

Since $C^{\infty}(\Omega)$ is dense both in $X$ and ${H_{\mu}(\Omega)}$, it follows that $X={H_{\mu}(\Omega)}$. Henceforth
we may suppose that the norm in $X$ coincides with the norm in
${H_{\mu}(\Omega)}$ and that the inner product in $X$ is given by
$<u,v>_X = <u,v>_{{H_{\mu}(\Omega)}}$. Let us note that the identification principle \cite[Identification Principle 21.18, pg. 254]{zei85}) implies that
if $<\cdot,\cdot>_{X,X^*}$ denotes the duality pairing on
$X$, then $<\cdot,\cdot>_{X,X^*}=<\cdot,\cdot> $.
To proceed further we note that
the bilinear form
\begin{eqnarray*}
\mathbf{a}(u,v) = \int_{\Omega} u\, v\,dx,\;\; \mbox{for all}\;\;
u,v \in X,
\end{eqnarray*}
is clearly continuous in $X$. The Riesz representation
theorem implies that we can define a bounded linear operator $\mathbf{L}$ such
that
\begin{eqnarray} \label{Rabi1.1}
\mathbf{a}(u,v) = <\mathbf{L}u,v>,\;\; \mbox{for all}\;\; u,v \in
X.
\end{eqnarray}
The operator $\mathbf{L}$ is self adjoint and compact and its largest
eigenvalue $\nu_1$ is characterized by
\begin{eqnarray*}
\nu_1= \sup_{u\in X} \frac{<\mathbf{L}u,u>}{<u,u>} = \sup_{u\in X}
\frac{\int_{\Omega} u^2\, dx}{\int_{\Omega} |\nabla u|^2\, dx -
\mu\, \int_{\Omega} \frac{u^2}{|x|^2}\, dx}.
\end{eqnarray*}
Then, by Lemma \ref{princeig} it readily follows that the positive
eigenfunction $u_1$ of (\ref{eq1.2}) corresponding to
$\lambda_{1,\mu}$ is a positive eigenfunction of $\mathbf{L}$
corresponding to $\nu_1=1/\lambda_{1,\mu}$.

After these preparations, we may define the
nonlinear operator $\mathbf{N}(\lambda,\cdot):\mathbb{R}\times
X\rightarrow X^*$ as
\begin{eqnarray}
\label{Rabi2} <\mathbf{N}(\lambda,u),v> = \int_{\Omega} \nabla u
\nabla v\, dx - \mu\, \int_{\Omega} \frac{u\, v}{|x|^2}\, dx
- \lambda\, \int_{\Omega} u\, v\, dx + \int_{\Omega}
|u|^{2\gamma}\, u\, v\,dx,
\end{eqnarray}
for all $v\in X$. Since the functional $S: X \to \mathbb{R}$
defined by
\begin{eqnarray*}
S(v)= \int_{\Omega} \nabla u  \nabla v\, dx - \mu\,
\int_{\Omega} \frac{u\, v}{|x|^2}\, dx - \lambda\, \int_{\Omega}
u\, v\, dx + \int_{\Omega} |u|^{2\gamma}\, u\, v\,dx,\; v\in X,
\end{eqnarray*}
is a bounded linear functional we have that $\mathbf{N}(\lambda,
u)$ is well defined from (\ref{Rabi2}). Moreover  by using the
fact that $X=H_{\mu}(\Omega)$ and relation (\ref{Rabi1.1}), we can
rewrite $\mathbf{N}(\lambda,u)$ in the form $\mathbf{N}(\lambda,u)
= u - \mathbf{G}(\lambda,u)$ where
$\mathbf{G}(\lambda,u):=\lambda\mathbf{L}u-\mathbf{H}(u)$,
\begin{eqnarray*}
<\mathbf{H}(u),v>=\int_{\Omega}|u|^{2\gamma}uv\,dx\;\;\mbox{for
all}\;\;v\in X.
\end{eqnarray*}
Under condition (\ref{cruc1}), the embedding $H^1_{\mu} (\Omega) \hookrightarrow
L^{2\gamma+2} (\Omega)$ is compact, implying that the map $\mathbf{H}$ is compact.
To check condition (\ref{Order}) of Theorem \ref{rab}, we derive first the inequality
\begin{eqnarray}
\label{limbif} \frac{1}{||u||_{X}}|<\mathbf{H}(u),v> |&\leq&
\frac{1}{||u||_{X}}||u||^{2\gamma}_{L^{2\gamma+2}}
||u||_{L^{2\gamma+2}}
||v||_{L^{2\gamma+2}}\nonumber\\
&\leq& c_1\, ||u||^{2\gamma}_{X} ||v||_{X}.
\end{eqnarray}
Then, we get from (\ref{limbif}) that
\begin{eqnarray*}
\lim_{||u||_{X}\rightarrow
0}\frac{||\mathbf{H}(u)||_{X^*}}{||u||_{X}}
=\lim_{{||u||_{X}\rightarrow
0}}\sup_{||v||_{X}\leq1}\frac{1}{||u||_{X}}|<\mathbf{H}(u),v>|=0.
\end{eqnarray*}

It remains to prove that  $C_{\lambda_{1,\mu}}$\ is global. We proceed in two steps \vspace{0.2cm}, adapting the arguments of \cite{kz06}. \\
(a)\ We shall prove first that all solutions\ $(\lambda,u) \in
C_{\lambda_{1,\mu}}$\ close to\ $(\lambda_{1,\mu} ,0)$ are positive for all
$x \in \Omega$. More precisely we shall prove that
there exists\ ${\epsilon}_0>0$,\ such that any
$(\lambda, u(x)) \in C_{\lambda_{1,\mu}} \cap B_{\epsilon_0}
((\lambda_{1,\mu},0))$,  satisfies $u(x) > 0$,\ for any\ $x \in
\Omega$. Here $B_{\epsilon_0} ((\lambda_{1,\mu},0))$, stands for the
open ball of $C_{\lambda_{1,\mu}}$ of center $(\lambda_{1,\mu},0)$
and radius $\epsilon_0$)

We argue by contradiction, assuming that  $(\lambda_n, u_n)$ is a sequence of
solutions of\ (\ref{eq1.1}), such that $(\lambda_n, u_n) \to
(\lambda_{1,\mu},0)$\ and that $u_n$\ are changing sign in
$\Omega$. Let $u_{n}^{-} := \min \{0, u_n\}$ and
${\mathcal{U}}_n^{-} =: \{ x \in \Omega: u_n(x)<0 \}$.\ Since
$u_n=u_n^+ -u_n^{-}$ is a solution of the problem (\ref{eq1.1}) it
can be easily seen that $u_n^{-}$, satisfies (in the weak sense)
the equation
\begin{eqnarray}
\label{Prtonos} - \Delta u_n^{-} - \mu\, \frac{u^{-}}{|x|^2} &-&
\lambda_n u_n^{-} + |u_n|^{2\gamma} u_n^{-} =0, \\
u_n^{-}|_{\partial \Omega} &=& 0. \nonumber
\end{eqnarray}
Then, multiplying $(\ref{Prtonos})$ with $u_{n}^{-}$ and
integrating over $\Omega$ we have that
\begin{eqnarray}
\label{Prtonos2} \int_{{\mathcal U}_n^{-}} | \nabla u_{n}^{-}
|^{2}\, dx  - \mu \int_{{\mathcal U}_n^{-}}
\frac{|u_{n}^-|^2}{|x|^2}\, dx - \lambda_n \int_{{\mathcal
U}_n^{-}} |u_{n}^-|^2\, dx + \int_{{\mathcal U}_n^{-}}
|u_n|^{2\gamma}|u_n^{-}|^2\, dx =0.
\end{eqnarray}
Since $\lambda_n$ is a bounded sequence, it follows  from
(\ref{Prtonos2}) and H\"{o}lder's inequality that
\begin{eqnarray}
\label{Prtonos3}
||u_{n}^{-}||^{2}_{H_{\mu}({\mathcal U}_n^{-})}  &\leq&
\lambda_n\, \int_{{\mathcal U}_n^{-}} |u_{n}^-|^2\, dx\nonumber\\
&\leq& C\,|{\mathcal U}_n^{-}|^{\frac{qN-2N+2q}{qN}}\,
\left(\int_{{\mathcal U}_n^{-}} |u_n|^{p^*}\right)^{\frac{2}{p^*}}\\
&\leq& C\, |{\mathcal U}_n^{-}|^{\frac{qN-2N+2q}{qN}}|
|u_{n}^-||^{2}_{H_{\mu}({\mathcal U}_n^{-})}.\nonumber
\end{eqnarray}
where $p*$ is the critical exponent defined in (\ref{KZ1}), for any $q\in [1,2)$. Then, from (\ref{Prtonos3}) we get that
\begin{equation} \label{sob}
M \leq |{\mathcal U}_n^{-}|,\;\;\mbox{for all}\;\; n,
\end{equation}
with the constant $M$ being independent of \ $n$. We denote now
by\ $\tilde{u}_n,=u_n/||u_n||$\ the normalization of\ $u_n$.\ Then
there exists a subsequence of\ $\tilde{u}_n$ (not relabelled)
converging weakly in\ $H_{\mu}(\Omega)$ \ to some function\
$\tilde{u}_0$. It can be seen that $\tilde{u}_0=u_{1,\mu}$.
Moreover, $\tilde{u}_n\to u_{1,\mu} >0$ in $L^2(\Omega)$. Passing
to a further subsequence if necessary,  by Egorov's Theorem,
$\tilde{u}_n \to u_{1,\mu} $ uniformly on $\Omega$ with the
exception of a set of arbitrary small measure. This contradicts
(\ref{sob}) and we conclude the functions \ $u_n$ cannot change
sign. \vspace{0.2cm} \\
(b)\ We shall exclude next that for some solution\ $(\lambda,u) \in
C_{\lambda_{1,\mu}}$, there exists a point\ $\xi \in \Omega$,\
such that\ $u(\xi) < 0$:
Using (a), the fact that the continuum
$C_{\lambda_{1,\mu}}$\ is connected (see Theorem \ref{rab}) and
the\ $C^{1,\zeta}_{loc}(\Omega \backslash \{0\})$- regularity of
solutions, we deduce that there exists\ $({\lambda}_0,u_0) \in
C_{\lambda_{1,\mu}}$,\ such that\ $u_0(x) \geq 0$,\ for all\ $x
\in \Omega$,\ except possibly some point\ $x_0 \in \Omega$,\
such that\ $u_0(x_0) = 0$.\ Then, the maximum principle (see
\cite{breduptes05, davdup03}) and the fact that the solutions are
singular at the origin imply that\ $u_0 \equiv 0$\ on\ $\Omega$.
Thus, we may construct a sequence\ $\{ ({\lambda}_n,u_n) \}
\subseteq C_{\lambda_{1,\mu}}$,\ such that\ $u_n(x)> 0$,\ for all\
$n$\ and\ $x \in \Omega$,\ $u_n \rightarrow 0$\ in\ $H_{\mu}
(\Omega)$,\ and\ ${\lambda}_n \rightarrow {\lambda}_0$.\ However,
this is true only for\ $\lambda_0 = \lambda_{1,\mu}$.\ As a
consequence, we have that\ $C_{\lambda_{1,\mu}}$\ cannot cross
$(\lambda,0)$ for some $\lambda \ne \lambda_1$, and every function
which belongs to $C_{\lambda_{1,\mu}}$ is strictly positive.
$\blacksquare$ \vspace{0.2cm}
\subsection{Approximation by bounded domains not containing the origin}
\label{THEOREMA2}
In this subsection we prove Theorem \ref{approx}. The proof is also
an alternative approach, to show Theorem \ref{globalbif}, approximating (\ref{eq1.1}) by the family of problems,
\[
(A)_r\;\;\;\;\;\; \biggl\{
\begin{array}{ll}
-\Delta u - \mu\, \frac{u}{|x|^2} = \lambda\, u - |u|^{2 \gamma}\,
u,\;\;\;\mbox{in}\;\;\Omega_r=\Omega \backslash B_r(0),
\vspace{0.2cm} \\ \;\;\;\;\;\;\;\;\;\;\;\;\;\; u |_{\partial
\Omega_r} = 0,
\end{array}
\]
for some $r>0$ sufficiently small. Standard regularity results
imply that if\ $u$\ is a weak solution of the problem ($(A)_r$),
for some $r>0$ small enough, then $u$ belongs at least in
$C^{1,\zeta}_{loc} (\Omega_r)$, for some\ $ \zeta \in (0,1)$.\
\vspace{0.1cm}

The corresponding approximating linear eigenvalue problems
\[
(AL)_r\;\;\;\;\;\; \biggl\{
\begin{array}{ll}
&-\Delta u - \mu\, \frac{u}{|x|^2} = \lambda\,
u,\;\;\;\mbox{in}\;\;\Omega_r, \vspace{0.2cm} \\
&\;\;\;u |_{\partial \Omega_r} = 0,
\end{array}
\]
admit for any $r>0$, a positive principal eigenvalue $\lambda_{1,\mu,r}$, characterized by
\[
\lambda_{1,\mu,r} = \inf_{
\begin{array}{c}
               \phi \in H_{0}^{1}(\Omega_r) \\
                \phi \not\equiv 0
               \end{array}}
\frac{\int_{\Omega_r} |\nabla \phi|^2\; dx - \mu\, \int_{\Omega_r}
\frac{|\phi|^2}{|x|^2}}{\int_{\Omega_r} |\phi|^2\; dx}.
\]
with the following properties: $\lambda_{1,\mu,r}$, is simple with
a positive associated eigenfunction $u_{1,\mu,r}$ and
$\lambda_{1,\mu,r}$ is the only eigenvalue of\ $(PL)_r$, with
positive associated eigenfunction. Furthermore, we have the following
\setcounter{lemma}{5}
\begin{lemma} \label{reigen}
Let $0<\mu\leq\mu^*$, and  $\lambda_{1,\mu}$ and $\lambda_{1,\mu,r}$, be the positive
principal eigenvalues of the problems (\ref{eq1.2}) and $(AL)_r$,
respectively. Then

(i)\ $u_{1,\mu,r}(x) \leq u_{1,\mu}(x),\;\;\; \mbox{for any}\; x
\in \bar{\Omega}_r,\;\; \mbox{and any}\; r>0$.

(ii) $u_{1,\mu,r} \to u_{1,\mu}$ in $H_{\mu}(\Omega) \cap
L^{\infty}_{loc} (\Omega\setminus\{0\})$, and\ $\lambda_{1,\mu,r}
\downarrow \lambda_{1,\mu}$,\ as\ $r \downarrow 0$.\
\end{lemma}
{\bf Proof:}\ (i) Having in mind, that both $u_{\lambda,r}$\ and
$u_{\lambda}$ are sufficiently smooth and positive functions on
$\bar{\Omega}_r$, the assertion follows
from the comparison principle (cf.
\cite[Theorem 10.5]{pucc04}).
\newline
(ii) We extend  $u_{1,\mu,r}$ on $\Omega$ as
\[
\hat{u}_{1,r} (x) = : \left\{
\begin{array}{ll}
u_{1,\mu,r} (x),\;\;\;& x \in \Omega_r, \\
0,\;\;\;& x \in B_r,
\end{array}
\right.
\]
for any sufficiently small\ $r>0$, using in the sequel for
convenience, the same notation
$u_{1,\mu,r}\equiv\hat{u}_{1,\mu,r}$. We note first that
\[
\lambda_{1,\mu,r} = \frac{\int_{\Omega_r} |\nabla u_{1,\mu,r}|^2\;
dx- \mu\, \int_{\Omega_r}
\frac{|u_{1,\mu,r}|^2}{|x|^2}}{\int_{\Omega_r} |u_{1,\mu,r}|^2\;
dx} = \frac{\int_{\Omega} |\nabla u_{1,\mu,r}|^2\; dx- \mu\,
\int_{\Omega} \frac{|u_{1,\mu,r}|^2}{|x|^2}}{\int_{\Omega}
|u_{1,\mu,r}|^2\; dx} \geq \lambda_{1,\mu}.
\]
Since $\Omega_{r_1} \subset \Omega_{r_2}$, for any $r_1 >
r_2$, we deduce that $\lambda_{1,\mu,r}$ is an decreasing sequence as $r \to 0$.
Moreover, $u_{1,\mu,r}$ forms  a bounded sequence in
$H_{\mu}(\Omega)$, thus $u_{1,\mu,r} \rightharpoonup u^*$ in
$H_{\mu}(\Omega)$ (up to a subsequence), and $\lambda_{1,\mu,r} \to \lambda^*$ in
$\mathbb{R}$ . Then, by the compact embedding $H_{\mu}
(\Omega) \hookrightarrow L^{2}(\Omega)$ we get that
\[
\lambda_{1,r} \int_{\Omega} |u_{1,r}|^2dx \to \lambda^*
\int_{\Omega} |u^*|^2dx,
\]
as\ $r \to 0$.\ Therefore,
\[
||u_{1,\mu,r}||_{H_{\mu}(\Omega)} \to ||u^*||_{H_{\mu}(\Omega)}.
\]
Hence $(\lambda^*,u^*)$ must be an
eigenpair of (\ref{eq1.2}) and from Lemma \ref{princeig} (ii),
we infer that \ $(\lambda^*,u^*) \equiv (\lambda_1,u_1)$.
Finally, we consider the difference $\psi = u - u_{\lambda,r}$.  Standard regularity
results imply that
\[
||\psi||_{W^{2,2}_{loc} (\Omega\setminus\{0\})} \leq C\,
||\psi||_{W^{1,2} (\Omega)} + O(r),\;\;\; \mbox{as}\;\; r \to 0,
\]
for some positive constant\ $C$\ independent from\ $r$. By
a bootstrap argument, we conclude that $u_{1,\mu,r} \to
u_{1,\mu}$ in\ $L^{\infty}_{loc} (\Omega\setminus\{0\})$.
$\blacksquare$
\vspace{0.2cm}

Rabinowitz's Theorem \ref{rab}, is applicable for the approximating problems $(A)_r$, by following closely the arguments used in proof of Theorem \ref{globalbif}.
\setcounter{lemma}{6}
\begin{lemma} \label{rglobalbif}
Assume that $0<\mu\leq\mu^*$. The principal eigenvalue $\lambda_{1,\mu,r}$\ of $(PL)_r$ is a
bifurcating point of the problem $(P)_r$ (in the sense of
Rabinowitz) and $C_{\lambda_{1,\mu,r}}$ is a global branch of
nonnegative solutions , which "bends" to the right of
$\lambda_{1,\mu}$. For any fixed $\lambda > \lambda_{1,\mu}$ these
solutions are unique.
\end{lemma}
The properties of the global branch $C_{\lambda_{1,\mu,r}}$ can be proved as in Proposition \ref{propbend} (see Subsection 2.4). The nonlinear analogue of Lemma \ref{reigen} is stated in
\setcounter{proposition}{7}
\begin{proposition} \label{rcomparison}
Assume that $0<\mu\leq\mu^*$, and let $\lambda$ be a fixed number, such that $(\lambda,
u_{\lambda,r}) \in C_{\lambda_{1,\mu,r}}$. Then,

(i)\ \ \ $u_{\lambda,r} \to u_{\lambda}$ in $H_{\mu}(\Omega)$,
with $(\lambda, u_{\lambda}) \in C_{\lambda_{1,\mu}}$,

(ii) $u_{\lambda,r}(x) \leq u_{\lambda}(x),\;\;\; \mbox{for any}\;
x \in \bar{\Omega}_r,\;\; \mbox{and any}\; r \downarrow 0$,

(iii) $u_{\lambda,r} \to u_{\lambda}$ in $L^{\infty}_{loc}
(\Omega\setminus\{0\})$, as $r \downarrow 0$.
\end{proposition}
{\bf Proof:}\ (i)\ We shall prove first that $u_{\lambda,r}$ is a
bounded sequence in $H_{\mu} (\Omega)$. We argue by contradiction, assuming that
\begin{eqnarray}
\label{eq3.8cont}
||u_{\lambda,r}||_{H_{\mu} (\Omega)} \to \infty\;\;\mbox{as}\;\;r
\downarrow 0.
\end{eqnarray}
From the weak formulation of the problems $(A)_r$ we get that $u_{\lambda,r}$ satisfies the equation
\begin{equation} \label{eq3.8}
\int_{\Omega_r} |\nabla u_{\lambda,r}|^2\; dx- \mu\,
\int_{\Omega_r} \frac{|u_{\lambda,r}|^2}{|x|^2} = \lambda\,
\int_{\Omega_r} |u_{\lambda,r}|^2\; dx - \int_{\Omega_r}
|u_{\lambda,r}|^{2\gamma+2}\; dx,
\end{equation}
which implies that
\begin{equation} \label{eq3.9}
||u_{\lambda,r}||_{H_{\mu} (\Omega)} \leq \lambda\,
||u_{\lambda,r}||_{L^2 (\Omega)},
\end{equation}
for any $r$ small enough. Setting
\[
\tilde{u}_{\lambda,r} =
\frac{u_{\lambda,r}}{||u_{\lambda,r}||_{H_{\mu} (\Omega)}},
\]
we get that $||\tilde{u}_{\lambda,r}||_{H_{\mu} (\Omega)} =1$,
for any $r>0$ small enough. Consequently (up to a subsequence) $\tilde{u}_{\lambda,r}$ converges
weakly to some $\tilde{u}_*$ in $H_{\mu} (\Omega)$, as $r \downarrow 0$,
and so $\tilde{u}_{\lambda,r} \to \tilde{u}_*$ in $L^2
(\Omega)$ as well as in $L^{2\gamma+2} (\Omega)$, as $r\downarrow 0$. In
addition, it follows from (\ref{eq3.9}) that
\[
||\tilde{u}_{\lambda,r}||_{H_{\mu} (\Omega)} \leq \lambda\,
||\tilde{u}_{\lambda,r}||_{L^2 (\Omega)},\;\;\;\; \mbox{for
any}\;\; r>0,
\]
hence $\tilde{u}_* \not\equiv 0$. Dividing (\ref{eq3.8}) by
$||u_{\lambda,r}||^{2\gamma+2}_{H_{\mu} (\Omega)}$ we get the equation
that
\begin{eqnarray}
\label{eq3.9b}
\int_{\Omega_r}
|\tilde{u}_{\lambda,r}|^{2\gamma+2}\; dx
=
\frac{\lambda}{||u_{\lambda,r}||^{2\gamma}_{H_{\mu} (\Omega)}}\,
\int_{\Omega_r} |\tilde{u}_{\lambda,r}|^2\; dx
-\frac{1}{||u_{\lambda,r}||^{2\gamma}_{H_{\mu} (\Omega)}} \left (
\int_{\Omega_r} |\nabla \tilde{u}_{\lambda,r}|^2\; dx- \mu\,
\int_{\Omega_r} \frac{|\tilde{u}_{\lambda,r}|^2}{|x|^2} \right ).
\end{eqnarray}
Passing to the limit to (\ref{eq3.9b}) as $r\downarrow 0$, (\ref{eq3.8cont}) implies that $\tilde{u}_{\lambda,r} \to 0$ in $L^{2\gamma+2}
(\Omega)$, contradicting that $\tilde{u}_*\not\equiv 0$.

Therefore, $u_{\lambda,r}$ is a bounded sequence in $H_{\mu} (\Omega)$
converging weakly to some $u_*$ in $H_{\mu}(\Omega)$
as $r \downarrow 0$. Then (\ref{eq3.9}), implies again that and $u_*\not\equiv 0$.
Passing to the limit to the weak formulation of (\ref{eq1.1}), we
deduce that $u_*$ is a
solution of (\ref{eq1.1}). We set $u_* = u_{\lambda}$.
Claims (ii) and (iii) can be proved by similar arguments to those used in Lemma
\ref{reigen}. $\blacksquare$ \vspace{0.2cm}

{\bf Proof of Theorem \ref{approx}:} We are making use of Whyburn's Theorem (see
\cite{eg00} and the references therein). For some $R>0$ and some
sequence $r_n \downarrow 0$, as $n \to \infty$, we define the sets
$A_n$ as follows:
\[
A_n = \biggr \{ B_R (\lambda_1 , 0) \cap
\mathcal{C}_{\lambda_{1,\mu,r_n}} \biggr \} \subset \mathbb{R}
\times H_{\mu} (\Omega).
\]
For every $n \in \mathbb{N}$, these sets are connected and closed.
In addition, Lemma \ref{rcomparison} implies that
\[
\liminf_{n \to \infty} \{A_n\}
\not\equiv \emptyset.
\]
We will justify  next, that the set $\bigcup_{n \in \mathbb{N}} A_n$
is relatively compact i.e., every sequence in $A_n$ contains a
convergent subsequence. To this end, we consider  $(\lambda_n,u_n) \in \bigcup_{n \in
\mathbb{N}} A_n$. Then, the sequence $(\lambda_n,u_n)$ is bounded
in $\mathbb{R} \times H_{\mu} (\Omega)$. Henceforth there exists a
subsequence still denoted by $(\lambda_n,u_n)$, such that
$\lambda_n \to \lambda_*$ and $u_n \rightharpoonup u_*$ in
$H_{\mu} (\Omega)$, $u_n\rightarrow u_*$ in $L^2(\Omega)$. Moreover,  $u_n$ satisfies (\ref{eq3.9}),
from which it readily follows that
$$||u_n||_{H_{\mu}(\Omega)}\rightarrow ||u^*||_{H_{\mu}(\Omega)}\;\;\mbox{as}\;\;n\rightarrow\infty.$$
Hence, the subsequence $u_n$ converges strongly to $u^*$ in $H_{\mu}(\Omega)$, and arguing as in Proposition \ref{rcomparison},  we get that $u_* \not\equiv 0$ as well as that $(\lambda_*,u_*)$ is a solution of
(\ref{eq1.1}). From the same token we have that
\[
\liminf_{n \to \infty} \{A_n\}=\limsup_{n \to \infty} \{A_n\}
\not\equiv \emptyset.
\]
Applying similar arguments, we may let $R \to \infty$ in order to obtain
that $\mathcal{C}_{\lambda_{1,\mu,r_n}}  \to C_{\lambda_1,\mu}$, in $\mathbb{R} \times
H_{\mu}(\Omega)$ for any $R \in
\mathbb{R}^+$. $\blacksquare$
\subsection{Behavior of the branch $C_{\lambda_{1,\mu}}$ as $\mu\rightarrow\mu^*$}
\label{THEOREMA3}
We conclude in this section, with the discussion on the  properties of the global branches
$C_{\lambda_{1,\mu}}$ when $0<\mu \leq \mu^*$. We start with the proof of Proposition \ref{propbend}, which actually shows that the global bifurcation is of supercritical type.
\vspace{0.2cm}
\newline
{\bf Proof of Proposition \ref{propbend}:}\ (i) Assume by contradiction that\ $C_{\lambda_{1,\mu}}$\ bends to
the left of\ $\lambda_{1,\mu}$. Then there exists a pair
$(\lambda,u) \in \mathbb{R} \times H_{\mu}(\Omega)$ with $0 < \lambda
< \lambda_{1,\mu}$, such that
\begin{equation} \label{eq2.1ab}
\int_{\Omega} |\nabla u|^2\, dx - \mu\, \int_{\Omega}
\frac{u^2}{|x|^2}\, dx = \lambda \int_{\Omega} |u|^2\, dx -
\int_{\Omega} |u|^{2\gamma+2}\, dx,
\end{equation}
Last equation implies that
\[
||u||^{2}_{H_{\mu}(\Omega)} \leq \lambda
||u||^{2}_{L^2(\Omega)},\;\;\; \mbox{with}\;\; \lambda<\lambda_1,
\]
contradicting the variational characterization of
$\lambda_{1,\mu}$. Thus, $C_{\lambda_{1,\mu}}$ must bend to the
right of $\lambda_{1,\mu}$. To show that $C_{\lambda_{1,\mu}}$ is
bounded for $\lambda$ bounded, we consider the weak formula satisfied by any $u\in C_{\lambda_{1,\mu}}$,
\begin{eqnarray}
\label{steadyA} \int_{\Omega}\nabla u\nabla\psi
dx-\int_{\Omega}\frac{u\psi}{|x|^2} dx-\lambda\int_{\Omega}u\psi
dx+\int_{\Omega}|u|^{2\gamma}u\psi dx=0,\;\;\mbox{for
all}\;\;\psi\in H_{\mu}(\Omega).
\end{eqnarray}
Setting $\psi=u$ in (\ref{steadyA}) and using the inequality
\begin{eqnarray}
\label{usualA}
2\lambda\int_{\Omega}|u|^{2}dx&\leq&
2\lambda|\Omega|^{\frac{\gamma}{\gamma
+1}}||u||^2_{L^{2\gamma+2}}
\leq\frac{1}{2}||u||^{2\gamma +2}_{L^{2\gamma
+2}}+R_0,\\ R_0&=&(2\lambda)^{\frac{\gamma
+1}{\gamma}}\frac{2^{\frac{1}{\gamma}}\gamma}{(\gamma +1)^{\frac{\gamma +1}{\gamma}}}|\Omega|\nonumber,
\end{eqnarray}
we get that any $u\in C_{\lambda_{1,\mu}}$, satisfies the  bound
\begin{eqnarray}
\label{usualB}
||u||^2_{H_{\mu}(\Omega)}\leq R_0.
\end{eqnarray}
The bound (\ref{usualB}), shows that any $u\in C_{\lambda_{1,\mu}}$,
is bounded for each fixed $\lambda$.

(ii)\ Let\ $u \in C_{\lambda_{1,\mu}}$,\ and suppose that\ $v$\ is
a nonnegative solution of (\ref{eq1.1}) with\ $u \not\equiv v$.
Considering the  approximating solutions $u_{\lambda,r}$ of the problems $(A)_r$, we get from Proposition
\ref{rcomparison} (ii) (comparison principle) that
\begin{eqnarray}
\label{bend2}
u_{\lambda,r}(x) \leq \min_{x \in \Omega} \{u(x),\; v(x)\}.
\end{eqnarray}
Then, by the $L^{\infty}_{loc}$-convergence of\ $u_{\lambda,r}$\
to\ $u$ of Lemma
\ref{rcomparison} (iii) and (\ref{bend2}), we infer that
\begin{eqnarray}
\label{bendB}
u(x) \leq v(x).
\end{eqnarray}
We apply next the weak formula (\ref{steadyA}) for the solutions $u$ and $v$, setting $\psi=v$ and $\psi=u$ respectively. Subtracting the resulting equations, we get that
\[
\int_{\Omega} (|u|^{2\gamma}v - |v|^{2\gamma}u)\, dx = 0,
\]
contradicting (\ref{bendB}), unless\ $u \equiv v$.\ $\blacksquare$
\vspace{0.2cm}

Finally, we discuss the behavior of the branches
$C_{\lambda_{1,\mu}}$, as $\mu \uparrow \mu^*$.
The eigenfunction
$u_{1,\mu^*}$ does not belong in $H_0^1 (\Omega)$, although the
eigenfunctions $u_{1,\mu}$, $0<\mu<\mu^*$, belong in $H_0^1
(\Omega)$. Therefore, the behavior of the branches
$C_{\lambda_{1,\mu}}$ as $\mu\uparrow\mu^*$ should be completely different if considered in $H_{\mu^*}(\Omega)$ and in $H_0^1(\Omega)$ respectively.
\vspace{0.2cm} \\
{\bf Proof of Theorem \ref{convergth}:} A.  By assumption, the pair  $(\lambda_n,u_n)$, satisfies
\begin{equation} \label{eq3.16}
\int_{\Omega} |\nabla u_n|^2\, dx - \mu_n\, \int_{\Omega}
\frac{|u_n|^2}{|x|^2}\, dx = \lambda_n\, \int_{\Omega} |u_n|^2\,
dx - \int_{\Omega} |u_n|^{2\gamma+2}\, dx,
\end{equation}
which implies that
\begin{equation} \label{eq3.16a}
\int_{\Omega} |\nabla u_n|^2\, dx - \mu_n\, \int_{\Omega}
\frac{|u_n|^2}{|x|^2}\, dx \leq \lambda_n\, \int_{\Omega}
|u_n|^2\, dx.
\end{equation}
On the other hand, by the definition of the $H_{\mu^*}(\Omega)$-norm and the hypothesis $\mu_n\uparrow\mu^*$, it follows that
\begin{eqnarray}
\label{eq3.16aa}
||u_n||_{H_{\mu^*} (\Omega)} \leq \int_{\Omega} |\nabla u_n|^2\, dx
- \mu_n\, \int_{\Omega} \frac{|u_n|^2}{|x|^2}\, dx.
\end{eqnarray}
Combining (\ref{eq3.16a}) and (\ref{eq3.16aa}) with the assumption that $|\lambda_n|\leq L$, we get the estimate
\begin{equation}\label{eq3.17}
||u_n||_{H_{\mu^*} (\Omega)} \leq \lambda_n\, ||u_n||^2_{L^2
(\Omega)} < L\, ||u_n||^2_{L^2 (\Omega)}.
\end{equation}
We employ an argument similar to the one used in the proof of Proposition \ref{rcomparison}, assuming by contradiction that $||u_n||_{H_{\mu^*} (\Omega)} \to \infty$ as $n \to
\infty$. We consider the normalization  $\hat{u}_n$  of
$u_n$ in $H_{\mu^*} (\Omega)$,
\[
\hat{u}_n = \frac{u_n}{||u_n||_{H_{\mu^*} (\Omega)}},
\]
which is a bounded sequence in $H_{\mu^*}(\Omega)$. Hence,  we may extract a subsequence (not relabelled),  converging weakly to some $\hat{u}_*$ in $H_{\mu^*} (\Omega)$. The compact embedding $H_{\mu^*} (\Omega)
\hookrightarrow L^2 (\Omega)$ and inequality (\ref{eq3.17}) imply
that $\hat{u}_* \not\equiv 0$. Dividing (\ref{eq3.16a}) by
$||u_n||_{H_{\mu^*} (\Omega)}^{2}$, we get the inequality
\begin{equation}\label{eq3.19}
\int_{\Omega} |\nabla \hat{u}_n|^2\, dx - \mu_n\, \int_{\Omega}
\frac{|\hat{u}_n|^2}{|x|^2}\, dx \leq \lambda_n\, \int_{\Omega}
|\hat{u}_n|^2\, dx < \infty.
\end{equation}
Moreover, dividing (\ref{eq3.16}) by
$||u_n||_{H_{\mu^*}(\Omega)}^{2\gamma+2}$, we get the equation
\begin{eqnarray}
\label{eq3.19b}
\int_{\Omega}
|\hat{u}_n|^{2\gamma+2}\, dx
=
\frac{\lambda_n}{||u_n||_{H_{\mu^*} (\Omega)}^{2\gamma}}\,
\int_{\Omega} |\hat{u}_n|^2\, dx - \frac{1}{||u_n||_{H_{\mu^*} (\Omega)}^{2\gamma}} \left (
\int_{\Omega} |\nabla \hat{u}_n|^2\, dx - \mu_n\, \int_{\Omega}
\frac{|\hat{u}_n|^2}{|x|^2}\, dx \right ).
\end{eqnarray}
Passing to the limit  to (\ref{eq3.19b}) as $n\rightarrow\infty$, we deduce that $u_* \equiv 0$, which is the
contradiction. Thus $u_n$ must be bounded in $H_{\mu^*} (\Omega)$,
and (up to some subsequence) converges weakly to some $u_*$ in
$H_{\mu^*} (\Omega)$.

The strong convergence $(\lambda_n, u_n)
\to (\lambda_*, u_*)$ in $\mathbb{R} \times H_{\mu^*} (\Omega)$, follows from the compactness of the embedding $H_{\mu^*} (\Omega)\hookrightarrow L^2 (\Omega)$ and (\ref{eq3.17}). Let us remark that if
$u_* \equiv 0$,  the same argument implies that $u_n \to 0$ in $H_{\mu^*} (\Omega)$. In this case, division of
(\ref{eq3.16}) by $||u_n||_{H_{\mu^*} (\Omega)}^{2}$ and passage to the limit, shows that $\lambda_n
\to \lambda_{1,\mu^*}$.

It remains to prove that the limit $(\lambda_*, u_*) \in
C_{\lambda_{1,\mu^*}}$. Note that for any $\phi \in C_0^{\infty} (\Omega)$,
\[
\int_{\Omega} \nabla u_n\nabla \phi\, dx - \mu^*\,
\int_{\Omega} \frac{u_n\, \phi}{|x|^2}\, dx - (\mu_n - \mu^*)\,
\int_{\Omega} \frac{u_n\, \phi}{|x|^2}\, dx = \lambda_n\,
\int_{\Omega} u_n\, \phi\, dx - \int_{\Omega} |u_n|^{2\gamma}\,
u_n\, \phi\, dx.
\]
Passing to the limit as $n\rightarrow\infty$, we need to show that the integral
\[
\int_{\Omega} \frac{u_n\, \phi}{|x|^2}\, dx,
\]
remains bounded for any $\phi \in C_0^{\infty} (\Omega)$ and any
$n\in\mathbb{N}$. This claim follows by H\"{o}lder's inequality and the continuous embedding $H_{\mu^*}(\Omega)\hookrightarrow L^{p^*}(\Omega)$, since
\begin{eqnarray}
\label{cruc1ain}
\left|\int_{\Omega} \frac{u_n\, \phi}{|x|^2}\, dx\right| \leq
||\phi||_{L^{\infty} (\Omega)}\, ||u_n||_{L^{p^*}(\Omega)}
\int_{\Omega} |x|^{-\frac{2Nq}{Nq-N+q}}\, dx.
\end{eqnarray}
The integral in the right hand side of (\ref{cruc1ain}) converges if $q>\frac{N}{N-1}$. Combining this requirement with (\ref{cruc1}), the condition (\ref{cruc1aN=3}) follows for the case $N=3$. When $N\geq 4$, the claim is valid under the condition (\ref{cruc1}). 
\vspace{0.2cm}\\
B. Let $(\lambda_n, u_n) \in C_{\lambda_{1,\mu_n}}$,
and assume that $\mu_n\uparrow\mu^*$ and $\lambda_n \to
\lambda_{1,\mu^*}$ as $n \to \infty$. Assuming further that $u_n$
remains bounded in $H_0^1 (\Omega)$, we may extract a subsequence still denoted by
$u_n$, which converges weakly to some $u_*$ in $H_0^1 (\Omega)$. Passing to the limit in the weak formula as $n\rightarrow\infty$, it
follows that $u_*, \mu^*, \lambda_{1,\mu^*}$, satisfy
\[
\int_{\Omega} \nabla u_*\nabla \phi\, dx - \mu^*\,
\int_{\Omega} \frac{u_*\, \phi}{|x|^2}\, dx = \lambda_{1,\mu^*}\,
\int_{\Omega} u_*\, \phi\, dx - \int_{\Omega} |u_*|^{2\gamma}\,
u_*\, \phi\, dx,
\]
for any $\phi \in C_0^{\infty} (\Omega)$. However, the variational
characterization of $\lambda_{1,\mu^*}$ implies that this is true
only if $u_* \equiv 0$. Therefore  $u_n \rightharpoonup 0$, in
$H_0^1 (\Omega)$. On the other hand, arguing as in part A., it can be seen from (\ref{eq3.16}) that the
normalization $\bar{u_n} = u_n / ||u_n||_{H_0^1 (\Omega)}$ converges (up to a subsequence)
weakly to $u_{1,\mu^*}$ in $H_0^1 (\Omega)$ which is
impossible. Thus, $u_n$ must be unbounded in $H_0^1 (\Omega)$. $\blacksquare$
\section{Definition of a gradient semiflow}
In this section we shall define a gradient semiflow associated to
the semilinear parabolic equation (\ref{eq1.0a}),
\begin{eqnarray}
\label{defsemiflow}
\mathcal{S}(t):H_{\mu}(\Omega)\rightarrow H_{\mu}(\Omega),\;\;0<\mu\leq\mu^*,
\end{eqnarray}
with
\begin{eqnarray}
\label{Liapunov}
\mathcal{J}(\phi) := \frac{1}{2}\, \int_{\Omega} |\nabla\phi|^2\,
dx - \frac{\mu}{2}\, \int_{\Omega}\frac{|\phi|^2}{|x|^2}\, dx -
\frac{\lambda}{2}\, \int_{\Omega} |\phi|^2\, dx +
\frac{1}{2\gamma+2}\, \int_{\Omega} |\phi|^{2\gamma +2}\, dx,\;\;0<\mu\leq\mu^*,
\end{eqnarray}
as a Lyapunov functional. In subsection \ref{SEMI2}, we discuss the stability properties of the equilibrium solutions by linearization. In subsection \ref{SEMI1}, and by following closely the general semiflow theory \cite{Ball2,jhale88,RTem88}, we  present the proof of Theorem \ref{main1}, as well as the description of the limit set $\omega(\phi_0)$ for nonnegative (nonpositive) initial data $\phi_0$, $\phi_0\not\equiv 0$ given in Corollary \ref{tsibadyoball}.
\subsection{Stability of equilibrium solutions by linearization}
\label{SEMI2}
\setcounter{equation}{0}
Seeking for nonpositive stationary solutions $u=-u_{-}$ with $u_{-}\geq 0$, $u_{-}\not\equiv 0$, it is clear that $u_{-}$ satisfies (\ref{eq1.1}). Therefore, Theorem \ref{globalbif}, can be restated as
\begin{corollary}
\label{nonposbif}
Let\ $\Omega \subset \mathbb{R}^N$,\ $N \geq 3$,\ be a bounded
domain. Assume that $0 < \mu\leq\mu^*$, and that condition (\ref{cruc1}) is satisfied.
Then, the principal eigenvalue $\lambda_{1,\mu}$\
of (\ref{eq1.2}) considered in $H^{1}_{\mu}(\Omega)$, is a bifurcating point of the problem
(\ref{eq1.1}) (in the sense of Rabinowitz) and
$C_{\lambda_{1,\mu}}$ and $C_{\lambda_{1,\mu}}^{-}$ are global branches of nonnegative and nonpositive  $H^{1}_{\mu}(\Omega)$- solutions respectively,
which bend to the right of $\lambda_{1,\mu}$. For any fixed
$\lambda > \lambda_{1,\mu}$, every solution $u \in C_{\lambda_{1,\mu}}$ and $u_{-} \in C_{\lambda_{1,\mu}}^{-}$ is the unique
nonnegative and unique nonpositive solutions for the problem (\ref{eq1.1}) and $u_{-}=-u$.
\end{corollary}

We first verify that solutions of (\ref{eq1.0a}) initiating from nonnegative (nonpositive) initial data remain nonnegative (nonpositive) for
all times. Then, we will proceed with the asymptotic stability of the nonnegative
equilibrium by linearization. For the latter, Hardy's inequalities and their improvements, allow for the definition of appropriate Garding forms, helping us to verify that zero is not an eigenvalue for the linearized flow around the nonnegative (nonpositive) equilibrium.
\setcounter{lemma}{1}
\begin{lemma}
\label{PositCone} Assume that $\mu \leq \mu^*$. The set
\[
\mathcal{D}_{+(-)}:=\left\{\phi\in
H_{\mu}(\Omega)\,:\,\phi(x)\geq (\leq ) 0\;\;\mbox{on}\;\;
\overline{\Omega}\right\},
\]
is a positively invariant set for the semiflow $\mathcal{S}(t)$.
\end{lemma}
{\bf Proof:}\ \ The argument of \cite[Proposition
5.3.1]{cazh} for the linear heat equation, can be repeated here (see also \cite{kz06}).  We assume that $\phi_0\in H_{\mu}(\Omega)$,
$\phi_0\geq 0$ a.e in $\Omega$, and
$\phi(t)=\mathcal{S}(t)\phi_0$, the global in time solution
of (\ref{eq1.0a}), initiating from $\phi_0$. We consider
$\phi^+:=\max\{\phi,0\}$, $\phi^{-}:=-\min\{\phi,0\}$. Both
$\phi^{+}$ and $\phi^{-}$ are nonnegative,  and
$\phi=\phi^+-\phi^{-}$.  It can be seen from (\ref{eq1.0a}) that
$\phi^{-}$ satisfies the equation
\begin{eqnarray} \label{eqminus}
\partial_t \phi^{-} - \Delta \phi^{-} - \mu\, \frac{\phi^{-}}{|x|^2}
- \lambda\, \phi^{-}+|\phi|^{2\gamma}\phi^{-}=0.
\end{eqnarray}
Moreover,  $\phi^{-}$ satisfies the energy equation (see Proposition \ref{main1}),
\begin{eqnarray}
\label{enegminus} \frac{1}{2}\frac{d}{dt}||\phi^{-}||^2_{L^2}+
\int_{\Omega} |\nabla\phi^{-}|^2\, dx - \mu\, \int_{\Omega}
\frac{|\phi^{-}|^2}{|x|^2}\, dx - \lambda\, ||\phi^{-}||^2_{L^2} +
\int_{\Omega} |\phi|^{2\gamma}|\phi^{-}|^2\, dx = 0.
\end{eqnarray}
From (\ref{enegminus}) and (\ref{varchareig}), we get that
\begin{eqnarray*}
\frac{1}{2}\frac{d}{dt}||\phi^{-}||^2_{L^2}\leq
c\,||\phi^{-}||^2_{L^2}.
\end{eqnarray*}
where $c=\lambda_{1,\mu}-\lambda$. Thus $\phi^{-}$ satisfies
\begin{eqnarray}
||\phi^{-}(t)||^2_{L^2}\leq
e^{ct}||\phi_0^{-}||^2_{L^2}=0,\;\;\mbox{for every}\;\;t\in
[0,+\infty),
\end{eqnarray}
implying that $\phi\geq 0$ for all $t\in (0,+\infty)$, a.e. in $\Omega$.
$\blacksquare$
\setcounter{proposition}{2}
\begin{proposition}
\label{Garding} Let $\mu \leq \mu^*$. The unique nonnegative (nonpositive)
equilibrium point which exists for $\lambda>\lambda_{1,\mu}$ is
uniformly asymptotically stable.
\end{proposition}
On the account of Corollary \ref{nonposbif}, we consider only the nonnegative equilibrium $u\geq 0$, $u\not\equiv 0$. First, we observe that the linearized semiflow around the zero solution, is defined by the Cauchy-Dirichlet problem
\begin{eqnarray*}
\label{linear0}
\partial_t \psi - \Delta u - \mu\, \frac{u}{|x|^2} &-&\lambda\,
\psi =0,\;\; x \in \Omega, \nonumber\\
\psi|_{\partial\Omega}&=&0.
\end{eqnarray*}
We have that $\phi=0$ is asymptotically stable in $H_{\mu}(\Omega)$
if $\lambda\leq\lambda_{1,\mu}$, and unstable in $H_{\mu}(\Omega)$
if $\lambda >\lambda_{1,\mu}$. The linearized semiflow around the nonnegative equilibrium point $u$ of
(\ref{eq1.0a}), is defined by the Cauchy-Dirichlet problem
\begin{eqnarray}
\label{linear}
- \Delta \psi - \mu\, \frac{\psi}{|x|^2} &-& \lambda\, \psi +
(2\gamma+1) |u|^{2\gamma} \psi=0,\\
\psi |_{\partial\Omega}&=&0,\nonumber
\end{eqnarray}
To confirm the asymptotic stability of $u$, we will prove that $\tilde{\mu}=0$, is not an eigenvalue for the eigenvalue problem
\begin{eqnarray}\label{Gard}
- \Delta \psi - \mu\, \frac{\psi}{|x|^2} &-& \lambda\, \psi +
(2\gamma+1) |u|^{2\gamma} \psi = \tilde{\mu} \psi,\\
\psi |_{\partial\Omega}&=&0.\nonumber
\end{eqnarray}
The weak formulation of (\ref{Gard}) is
\begin{eqnarray} \label{weakgard}
A_{\mu}(\psi,\omega) &:=& \int_{\Omega} \nabla \psi\nabla
\omega\, dx - \mu \int_{\Omega} \frac{\psi\, \omega}{|x|^2}\, dx
- \lambda\, \int_{\Omega} \psi\, \omega\, dx \nonumber \\
&+&(2\gamma+1)\, \int_{\Omega} |u|^{2\gamma} \psi\, \omega\,dx=
\tilde{\mu}\, \int_{\Omega} \psi\, \omega\, dx,
\end{eqnarray}
for every $\omega\in {H}_{\mu}(\Omega)$. Using the improved Hardy's inequality and the properties of the $H_{\mu}(\Omega)$-space for any $0<\mu\leq\mu^*$, we may consider a symmetric bilinear
form $A_{\mu} : {H}_{\mu}(\Omega) \times {H}_{\mu}(\Omega)
\rightarrow \mathbb{R}$, which in turns,  defines  a Garding form \cite[pg.
366]{zei85}:  Since
\begin{eqnarray*}
A_{\mu}(\psi,\psi)\geq ||\psi||^2_{{H}_{\mu}(\Omega)} - \lambda\,
||\psi||^2_{\mathrm{L^2(\Omega)}},
\end{eqnarray*}
Garding's inequality is satisfied. Then, it follows from
\cite[Theorem 22.G pg. 369-370]{zei85} and (\ref{HaPoiim4}), that the problem
(\ref{Gard}) has infinitely many eigenvalues of finite
multiplicity. Counting  the eigenvalues according to their
multiplicity, we derive the sequence
\begin{eqnarray}
- \lambda < \tilde{\mu}_1 \leq \tilde{\mu}_2 \leq \cdots,\;\;
\mbox{and}\;\; \tilde{\mu}_j \rightarrow \infty\;\; \mbox{as}\;\;
j \rightarrow \infty.
\end{eqnarray}
The smallest eigenvalue can be characterized by the minimization
problem
\begin{eqnarray} \label{minlin1}
\tilde{\mu}_1 = \min A_{\mu} (\psi,\psi),\;\; \psi \in H_{\mu}
(\Omega),\;\; ||\psi||_{L^2}=1.
\end{eqnarray}
The $j$-th eigenvalue, can be characterized by the minimum-maximum
principle
\begin{eqnarray} \label{minlin2}
\tilde{\mu}_j= \min_{M \in \mathcal{L}_j} \max_{\psi \in M}
A_{\mu}(\psi,\psi).
\end{eqnarray}
where $M=\{\psi\in H_{\mu} (\Omega)\, :\, ||\psi||_{L^2}=1 \}$ and
$\mathcal{L}_j$ denotes the class of all sets $M\cap L$ with $L$
an arbitrary $j$-dimensional linear subspace of $H_{\mu}(\Omega)$.

By using similar arguments as for the proof of Lemma
\ref{princeig} (see also Lemma \ref{reigen}), we may see that for
(\ref{Gard}), the (nontrivial) eigenfunction corresponding to the
principal eigenvalue $\tilde{\mu}_1$ is nonnegative, i.e
$\psi_1\geq 0$ a.e. on $\Omega$. Since $\tilde{\mu}_1,\psi_1$
satisfy (\ref{weakgard}) we get by setting $\omega=u$ that
\begin{eqnarray*}
\int_{\Omega} \nabla\psi_1\nabla u\, dx - \mu \int_{\Omega}
\frac{\psi_1\, u}{|x|^2}\, dx - \lambda\, \int_{\Omega} \psi_1\,
u\,dx + (2\gamma +1)\, \int_{\Omega} |u|^{2\gamma}\, \psi_1\, u\,
dx = \tilde{\mu}_1 \int_{\Omega} \psi_1\, u\, dx.
\end{eqnarray*}
On the other hand, by setting $\psi=\psi_1$ to the weak formula
(\ref{steadyA}) we get
\begin{eqnarray*}
\int_{\Omega} \nabla \psi_1\nabla u\, dx - \mu
\int_{\Omega} \frac{\psi_1\, u}{|x|^2}\, dx - \lambda\,
\int_{\Omega} \psi_1\, u\,dx + \int_{\Omega} |u|^{2\gamma}\,
\psi_1\, u\,dx=0.
\end{eqnarray*}
Subtracting these equations, we obtain that
\begin{eqnarray}
2\gamma\, \int_{\Omega} |u|^{2\gamma}\, u\, \psi_1\, dx =
\tilde{\mu}_1\, \int_{\Omega} u\, \psi_1\, dx,
\end{eqnarray}
which implies that $\tilde{\mu}_1>0$. Thus $\tilde{\mu}=0$ is not an eigenvalue, and $u$ is uniformly asymptotically stable.  $\blacksquare$
\subsection{Global attractor in $H_{\mu}(\Omega)$ for any $0<\mu\leq\mu^*$}
\label{SEMI1}
The proof of  Proposition \ref{main1} is based on the analogue of \cite[Theorem 3.6, pg. 40]{Ball2},  this time for the parabolic equation (\ref{eq1.0a}).
\vspace{0.2cm}
\newline
{\bf Proof of Proposition \ref{main1}:} It follows from \cite{vz00}, that the operator
$\mathcal{L}=-\Delta-\frac{\mu}{|x|^2}$ with domain of definition
(\ref{HaPoi2}) is a generator of a strongly continuous semigroup
$\mathcal{T}(t)$, for any $0<\mu\leq\mu^*$,  while the function $f(s)=|s|^{2\gamma}s-\lambda s$, defines a locally Lipschitz map $f: H_{\mu}(\Omega)\rightarrow L^2(\Omega)$ as it can be easily deduced by Lemma \ref{aux1}. Thus for any $0<\mu\leq\mu^*$ and any $\phi_0\in H_{\mu}(\Omega)$, there exists a unique solution $\phi$ of (\ref{eq1.0a}), defined on a maximal interval $[0,T_{max})$ and in the class
$\mathrm{C}([0,T];H_{\mu}(\Omega)\cap
\mathrm{C}^1([0,T];L^2(\Omega))$. The solution satisfies the variation
of constants formula
\begin{eqnarray}
\label{pmild}
\phi(t)=\mathcal{T}(t)\phi_0+\int_{0}^{t}\mathcal{T}(t-s)f(\phi(s))
ds.
\end{eqnarray}
Lemma \ref{aux1} implies also that the functional $\mathcal{J}\in C^1(\mathbb{R}, H_{\mu}(\Omega))$,  for any $0<\mu\leq\mu^*$.
Moreover, for all $\phi\in D(\mathcal{L})$ and any $t\in [0, T]$, $T<T_{max}$,
\begin{eqnarray}
\label{preliap} \left<\Delta\phi+\mu\frac{\phi}{|x|^2}+f(\phi),
\mathcal{J}'(\phi)\right>=-\int_{\Omega}\left|\Delta\phi+\mu\frac{\phi}{|x|^2}+f(\phi)\right|^2dx
=-\int_{\Omega}|\partial_t\phi|^2dx\leq 0.
\end{eqnarray}
Setting $h(t)=f(\phi(t))$, we consider the sequence $h_n(t)\in
C^{1}([0,T];H_{\mu}(\Omega))$ and $\phi_{0n}\in
D(\mathcal{L})$ such that
\begin{eqnarray*}
h_n&\rightarrow& h,\;\;\mbox{in}\;\;C^{1}([0,T];H_{\mu}(\Omega)),\\
\phi_{0n}&\rightarrow&\phi_0,\;\;\mbox{in}\;\; H_{\mu}(\Omega).
\end{eqnarray*}
We define
$\phi_n(t)=\mathcal{T}(t)\phi_{0n}+\int_{0}^{t}\mathcal{T}(t-s)h_n(s)ds$,
and it follows from \cite[Corrolary 2.5, p107]{Pazy83} that
$\phi_n(t)\in D(\mathcal{L})$, $\phi_n\in
C^{1}([0,T];H_{\mu}(\Omega))$ and that they satisfy
\begin{eqnarray}
\label{seqeq}
\partial_t\phi_n-\Delta\phi_n-\mu\frac{\phi_n}{|x|^2}+f(\phi_n)=0.
\end{eqnarray}
Moreover, from \cite[Lemma 5.5, pg. 246-247]{Ball1a} or
\cite[Theorem 3.6, pg. 41]{Ball2}, we deduce that
\begin{eqnarray*}
\phi_n\rightarrow\phi,\;\;\mbox{in}\;\;H_{\mu}(\Omega).
\end{eqnarray*}
Finally, by using the continuity of $\mathcal{J}$ and
(\ref{preliap}), and  passing to the limit to the equation
\begin{eqnarray*}
\mathcal{J}(\phi_n(t))-\mathcal{J}(\phi_{0n})&=&\int_{0}^{t}
\left<\mathcal{J}'(\phi_n(s)),\Delta\phi_n(s)+\mu\frac{\phi_n(s)}{|x|^2}+h_n(s)\right>ds \\
&=&-\int_{0}^{t}||\partial_t\phi_n(s)||^2_{L^2}ds+\int_{0}^{t}
\left<\mathcal{J}'(\phi_n(s)),h_n(s)-f(\phi_n(s))\right>ds,
\end{eqnarray*}
we derive
\begin{eqnarray}
\label{liap1}
\frac{d}{dt}\mathcal{J}(\phi(t))=-\int_{\Omega}|\partial_t\phi|^2dx,\;\;\mbox{for all}\;\;0<\mu\leq\mu^*\;\;\mbox{and}\;\;t\in [0,T],\;T<T_{max}.
\end{eqnarray}
From (\ref{liap1}) we infer that the unique solution $\phi$,
satisfies the {\em energy
equation}
\begin{eqnarray}
\label{eneg2} \frac{1}{2}\frac{d}{dt}||\phi||^2_{L^2}+
\int_{\Omega}|\nabla\phi|^2dx-\mu\int_{\Omega}\frac{|\phi|^2}{|x|^2}-\lambda||\phi||^2_{L^2}
+\int_{\Omega}|\phi|^{2\gamma +2}dx=0,\;\;\mbox{for all}\;\;0<\mu\leq\mu^*.
\end{eqnarray}
When $\lambda\leq\lambda_{1,\mu}$, we observe by using (\ref{varchareig}), that $\limsup_{t\rightarrow\infty}||\phi(t)||_{L^2}^2=0$. For the case $\lambda>\lambda_{1,\mu}$, we insert (\ref{usualA})
to (\ref{eneg2}), to get the estimate
\begin{eqnarray*}
\label{eneg3} \frac{1}{2}\frac{d}{dt}||\phi||^2_{L^2}+
\frac{1}{2}||\phi||_{H_{\mu}(\Omega)}^2+\lambda||\phi||^2_{L^2}
+\frac{1}{2}||\phi||^{2\gamma +2}_{L^{2\gamma +2}} dx\leq
\mathrm{R}_0.
\end{eqnarray*}
Then by Gronwall's Lemma
 \begin{eqnarray}
 \label{abset1}
 ||\phi(t)||^2_{L^2}\leq||\phi(0)||^2_{L^2}\exp(-2\lambda t)+
\frac{\mathrm{R}_0}{\lambda}(1-\exp(-2\lambda t)).
 \end{eqnarray}
Letting $t\rightarrow \infty$, from  (\ref{abset1}) we obtain that
\begin{eqnarray}
\label{abset2}
\limsup_{t\rightarrow\infty}||\phi(t)||^2_{L^2}\leq\rho^2,\;\;\;\rho^2=\mathrm{R}_0/\lambda.
\end{eqnarray}
Now assume that $\phi_0$ is in a bounded set $\mathcal{B}$ of
$H_{\mu}(\Omega)$. Then  (\ref{abset2}) implies that for any
$\rho_1>\rho$, there exists $t_0(\mathcal{B},\rho_1)$, such that
\begin{eqnarray}
\label{boundV}
||\phi(t)||^2_{L^2}\leq \rho_1^2,\;\;\mbox{for any}\;\;t\geq t_0(\mathcal{B},\rho_1).
\end{eqnarray}
By the definition of the Lyapunov functional $\mathcal{J}$ and (\ref{boundV}), we have the inequality
\begin{eqnarray}
\label{abset3} \mathcal{J}(\phi(t))&\geq&
\frac{1}{2}\int_{\Omega}|\nabla\phi|^2dx-
\frac{\mu}{2}\int_{\Omega}\frac{|\phi|^2}{|x|^2}dx-\frac{\lambda}{2}\int_{\Omega}|\phi|^2dx\nonumber\\
&\geq&
\frac{1}{2}\int_{\Omega}|\nabla\phi|^2dx-\frac{\mu}{2}\int_{\Omega}\frac{|\phi|^2}{|x|^2}dx
-\frac{\lambda}{2}\rho_1^2,\;\;t\geq
t_0.
\end{eqnarray}
Since $\mathcal{J}$ is
nonincreasing in $t$, we conclude with the bound
\begin{eqnarray}
\label{boundV1}
||\phi(t)||_{H_{\mu}(\Omega)}^2\leq 2\mathcal{J}(\phi_0)+\lambda\rho_1^2,\;\;t\geq t_0.
\end{eqnarray}
establishing that solutions are globally defined in
$H_{\mu}(\Omega)$, for any $0<\mu\leq\mu^*$ and $\lambda>\lambda_{1,\mu}$.
In addition, (\ref{boundV1}) implies that the semiflow $\mathcal{S}(t)$ is eventually bounded and since the operator
$\mathcal{L}$ has compact resolvent, $\mathcal{S}(t)$ is asymptotically compact (cf. \cite[Proposition 2.3, pg. 36]{Ball2}, \cite{jhale88,RTem88}). Thus, the positive orbit
$\gamma^+(\phi_0)$ for any $\phi_0\in H_{\mu}(\Omega)$ is
precompact and has a nonempty compact and connected invariant
$\omega$-limit set $\omega(\phi_0)$. Moreover (\ref{liap1}) implies that $\omega(\phi_0)\in\mathcal{E}$. Equilibria of $\mathcal{S}(t)$ are extreme points of $\mathcal{J}$, satisfying the weak formula (\ref{steadyA}).  From  (\ref{usualB}), we have that $\mathcal{E}$ is bounded for any fixed $\lambda$.  Hence $\mathcal{S}(t)$ is point dissipative.\ $\blacksquare$
\vspace{0.2cm} \\
{\bf Proof of Theorem \ref{tsibadyoball}:} Lemma \ref{PositCone}
and Proposition \ref{main1}, imply that the solution $\phi
(t)=\mathcal{S}(t)\phi_0$, initiating from initial data
$\phi_0\geq 0$ ($\phi_0\leq 0$), $\phi_0\not\equiv 0$ converge
towards the set of nonnegative (nonpositive) solutions of
(\ref{eq1.1}) as $t\rightarrow\infty$, in $H_{\mu} (\Omega)$, for
any $0<\mu\leq\mu^*$.  In fact, it follows from Theorem
\ref{globalbif} that the set of equilibrium solutions
$\mathcal{E}=\{u_{-},0, u\}$, when  $\lambda > \lambda_{1,\mu}$,
the trivial solution being unstable by Proposition \ref{Garding}.
Thus for any nonnegative (nonpositive) initial condition $\phi_0$,
$\omega(\phi_0)=\{u\}$ ($\omega(\phi_0)=\{u_{-}\}$). While in the
case $\lambda < \lambda_{1,\mu}$, Theorem \ref{globalbif} combined
with Propositions \ref{Garding} and \ref{main1}, imply that
$\mathrm{dist}(\mathcal{S}(t)\mathcal{B}, \{0\})\rightarrow 0$ as
$t\rightarrow\infty$, for every bounded set $\mathcal{B}\subset
H_{\mu} (\Omega)$. Thus, when $\lambda < \lambda_{1,\mu}$ the
global attractor $\mathcal{A}=\{0\}$. $\blacksquare$
\vspace{0.2cm}
\newline
{\bf Acknowledgements.} We would like to thank Prof. Anargyros Delis for stimulating discussions.
\bibliographystyle{amsplain}

\end{document}